\documentclass[12pt]{amsart}
\usepackage{amscd,amsfonts,amsmath,amssymb,amsthm,amstext,latexsym}
\usepackage{enumerate,pifont}
\usepackage[english]{babel}
\def\a{{\alpha}}

\def\g{{\gamma}}
\def\b{{\beta}}

\newcommand{\Z}{{\mathbb Z}}

\newcommand{\R}{{\mathbb R}}
\newcommand{\C}{{\mathbb C}}
\newcommand{\N}{{\mathbb N}}

\def\t{{\theta}}

\def\G{{\Gamma}}

\newtheorem{theorem}{Theorem}%[section]
\newtheorem{lemma}[theorem]{Lemma}
\newtheorem{proposition}[theorem]{Proposition}

\newtheorem{corollary}[theorem]{Corollary}

\begin{document}
\title{Index and nullity of the Gauss map of the Costa-Hoffman-Meeks surfaces}
\address{Universit\'e Paris-Est, 
Laboratoire d'Analyse et Math\'ematiques Appliqu\'ees,
5 blvd Descartes, 77454 Champs-sur-Marne, FRANCE}
\email{filippo.morabito@univ-mlv.fr}
\address{Universit\`a Roma Tre, Dipartimento di Matematica,  
Largo S. L. Murialdo 1, 00146 Roma, ITALIA}
\email{morabito@mat.uniroma3.it}
\author{Filippo Morabito}

\begin{abstract} 
The aim of this work is to extend the results
 of S. Nayatani about the index and the nullity of the Gauss map of the 
 Costa-Hoffman-Meeks surfaces for values of the genus bigger than 37.
 That allows us to state that these minimal surfaces are non degenerate for all
 the values of the genus in the sense of the definition of J. P\'erez and A. Ros.
 \end{abstract}
\maketitle

\section*{Introduction}
In the years 80's and 90's the study of the index of minimal 
surfaces in Euclidean space has been quite active. D. 
Fischer-Colbrie in \cite{F}, R. Gulliver and H. B. Lawson in 
\cite{GL} proved independently that a complete minimal 
surface $M$ in ${\mathbb R}^{3}$ with Gauss map $G$ has 
finite index if and only if it has finite total curvature. 
D. Fischer-Colbrie also observed that if $M$ has finite total 
curvature its index coincides with  the index of an operator 
$L_{\bar G}$ (that is the number of its
negative eigenvalues) associated to the extended Gauss map 
$\bar G$ of $\bar M,$ the compactification of $M.$ 
Moreover $N(\bar G),$ the null space of 
$L_{\bar G},$ if restricted to $M$ consists of the bounded 
solutions of the Jacobi equation. 
The nullity, ${\rm Nul}(\bar G),$ that is the 
dimension of $N(\bar G),$ and 
the index are invariants of $\bar G$ because they are 
independent of the choice of the conformal metric 
on $\bar M.$ \\

\noindent
The computation of the index and of the nullity of the Gauss map 
of the Costa surface and of the Costa-Hoffman-Meeks surface of
 genus $g=2,\ldots,37$ appeared respectively in the works \cite{N2} 
and \cite{N1} of S. Nayatani. The aim of this work is to extend his 
results to the case where $g \geqslant 38.$  \\

\noindent
In \cite{N2} he studied the index and the nullity of the 
 operator $L_G$ associated to an arbitrary holomorphic 
 map $G: \Sigma \to S^2,$
 where $\Sigma$ is a compact Riemann surface. 
He considered a deformation $G_t : \Sigma \to S^2,$ $t\in (0,+\infty),$
with $G_1=G$ (see equation \eqref{Gt}) and gave lower and upper bounds 
for the index of $G_t,$  
${\rm  Ind}(G_t),$ and its nullity, ${\rm Nul}(G_t),$ for $t$ near to $0$ 
and $+\infty$ 
and $t=1.$ 
The computation of the index and the nullity in the case of the Costa surface
is based on the fact that the Gauss map of this surface is a deformation for a particular 
value of $t$ of  the map $G$ defined by $\pi \circ G=1 / \wp',$ that is its stereographic 
projection is equal to the inverse of the derivative of the Weierstrass 
$\wp$-function for a unit 
square lattice. 
S. Nayatani computed ${\rm Ind}(G_t)$ and ${\rm Nul}(G_t)$
for  $t\in (0,+\infty),$ where $G$ is the map defined above. 
So the result concerning the Costa surface follows as a simple 
consequence from that. He obtained that for this surface the index 
and the  nullity are equal respectively to 5 and 4.\\

\noindent
In \cite{N1} S. Nayatani extended the last result  treating the case of 
the Costa-Hoffman-Meeks surface of genus $g,$ 
$M_g,$ %AGGIUNTO SIMBOLO PER DARE ENUNCIATI.
but only  for 
$2 \leqslant g \leqslant 37.$ 
He obtained that the index is equal to $2g+3$ and the nullity is equal to $4.$
Here we will show that these results continue to hold 
also for $g \geqslant 38.$ \\ %In fact we will prove the following theorems. \\

\noindent
J. P\'erez and A. Ros in \cite{PR} introduced a notion 
of non degenerate  minimal surface in terms of 
the  Jacobi functions having logarithmic growth at the ends 
of the surface.  
As consequence of the works \cite{N1} and \cite{N2},
the Costa-Hoffman-Meeks surface was known
to be non degenerate %with respect to the definition,
but only for $1 \leqslant g \leqslant 37.$\\

\noindent
The result of S. Nayatani about the nullity of the Gauss map of the 
Costa-Hoffman-Meeks surface  is essential 
for the construction due to L. Hauswirth and F. Pacard \cite{HP} of a 
family of  minimal surfaces  with two limit ends
asymptotic  to half Riemann minimal surfaces and of genus $g$ with 
$1\leqslant g \leqslant 37.$ Their construction is based on 
a gluing procedure
which involves the Costa-Hoffman-Meeks surface of genus $g$
and two half Riemann minimal surfaces. 
In particular the authors needed show the existence of a family 
of minimal surfaces
close to the  Costa-Hoffman-Meeks surface,  invariant under the 
action of the symmetry
with respect to the vertical plane $x_2=0,$ having one 
horizontal end asymptotic to 
the plane $x_3=0$ and having the upper 
and the  lower end asymptotic (up to translation) respectively 
to the upper and the lower
end of the standard catenoid whose axis of revolution is directed 
by the vector $\sin \t\, e_1+\cos \t\, e_3,$ $\t \leqslant \t_0$ 
with $\t_0$ sufficiently small.
That was obtained by Sch\"auder fixed point theorem 
and using the fact that the nullity of the Gauss map of the surface 
is equal to $4.$ In \cite{HP} the authors refer to this last result as 
a non degeneracy property of the Costa-Hoffman-Meeks surface.
It is necessary to remark that here the choice of working with 
symmetric deformations of the surface with respect to 
the plane  $x_2=0,$ has a key role.
Because of the restriction on the value of the genus which affects 
the result of S. Nayatani,
it was not possible to prove the existence of this family of minimal 
surfaces for higher values of the genus.\\

\noindent
So one of the consequences of our work is  the proof of
the non degeneracy of the Costa-Hoffman-Meeks surface for 
$g \geqslant 1$ in the sense of the definition given in \cite{PR} 
and also, only in a symmetric setting, in \cite{HP}.
So we can state that the family of examples constructed by
L. Hauswirth and F. Pacard exists for all the values of the genus.\\
%Moreover our result allows us to show in \cite{M} 
%the existence of a family of minimal 
%deformations of the  Costa-Hoffman-Meeks surface 
%for each value of the genus.\\

\noindent
Summarizing we will prove the following theorems.

\begin{theorem}
\label{indice.costa}
For $1 \leqslant g <+\infty$ 
the index of the Gauss map of $M_g$ is equal to $2g+3.$    
\end{theorem}

\begin{theorem}
\label{teorema.nullity}
For $1 \leqslant g <+\infty$ the null space of the Jacobi operator of $M_g$
has dimension equal to $4.$
\end{theorem}

\noindent
Using the definition of non degeneracy given in \cite{PR}, we can also rephrase 
this last result giving the following statement. 

\begin{corollary}
\label{corollario.nondeg}
The surface $M_g$ is non degenerate for  $1\leqslant g < \infty.$
\end{corollary}

\noindent
The author wishes to thank S. Nayatani for having provided the 
background computations on which are based the results about the 
Costa-Hoffman-Meeks surfaces contained in \cite{N1}.\\   

\noindent
The author is grateful to his thesis director, L. Hauswirth,
 for his support and for having brought this problem to his
attention.\\

\section{Preliminaries}
Let $M$ be a complete oriented minimal surface in ${\mathbb R}^3.$
The Jacobi operator of $M$ is $$L=-\Delta + 2 K$$ where $\Delta$ is the 
Laplace-Beltrami operator and $K$ is the Gauss curvature. 
Moreover we suppose that $M$ has finite total curvature.
Then $M$ is conformally equivalent to a compact Riemann surface with finitely 
many punctures and the Gauss map $G:M\to S^2$ extends to the compactified surface 
holomorphically. So in the following we will pay attention to a generic compact Riemann 
surface, denoted by  $\Sigma$ and $G:\Sigma \to S^2$ a non constant holomorphic map, where $S^2$ is the unit sphere in ${\mathbb R}^3$ endowed with the complex 
structure induced by the stereographic projection from the north pole (denoted by $\pi$). 
We fix a conformal metric $ds^2$ on $\Sigma$  and consider the operator
$L_G=-\Delta +|dG|^2, $ acting on functions on $\Sigma.$\\

\noindent
We denote by $N(G)$ the kernel of $L_G.$
We define ${\rm Nul}(G),$ the nullity of $G,$ as the dimension of $N(G).$ Since 
$L(G)=\{a\cdot G\,|\,a \in {\mathbb R}^3\}$ is a three dimensional subspace of 
 $N(G),$ then ${\rm Nul}(G) \geqslant 3.$ We denote the index
 of $G,$ that is the number of negative eigenvalues
 of $L_G,$ by ${\rm Ind}(G).$  The index and the nullity 
 are invariants of the map $G$: they are independent of
 the metric on the surface $\Sigma.$ So we can consider on $\Sigma$
 the metric induced by $G$ from $S^2.$ \\

\noindent
 N. Ejiri and M. Kotani in \cite{EK} and S. Montiel and A. Ros in \cite{MR}
  proved that 
a non linear element of $N(G)$ is expressed as the support function  of a 
complete branched minimal surface with planar ends whose extended Gauss map
is $G.$ In the following we will review briefly some results contained in \cite{MR}
used by S. Nayatani in \cite{N2}.\\

\noindent
We will use some  definitions and concepts
of algebraic geometry. They are recalled in 
subsection \ref{RiemannRoch}.\\  

 \noindent
 Let $\g$ be the meromorphic function defined by $\pi \circ G.$ 
 Let $p_j$ and $r_i$ be respectively the poles and the branch points of $\g.$ 
We denote by  $P(G)=\sum_{j=1 }^\nu n_j p_j,$  $S(G)=\sum_{i=1 }^\mu m_i r_i$ 
respectively the polar and ramification
 divisor of $\g.$ Here $n_j,m_i$ denote, respectively, the multiplicity 
of the pole $p_j$ and the multiplicity with which $\g$ takes 
its value at $r_i.$
  We define on the surface $\Sigma$ the divisor
 $$ D(G)=S(G)-2 P(G)$$
% and introduce the vector spaces $\hat H(G)$ and $H(G)$ by 
 and introduce the vector space $\bar H(G)$ (see \cite{MR}, theorem 4)

% \begin{equation}
% \label{Hcappello.G}
% \hat H(G)= \left\lbrace  \omega \in H^0(k_\Sigma+D(G)) \,| \, Res_{r_i} \omega=0, 
% 1 \leqslant i \leqslant \mu\right\rbrace,
% \end{equation}
% 
% $$H(G)=\left\lbrace  \omega \in \hat H(G) \, | \,Re \int_\a (1-g^2,i(1+g^2),2g)\o=0,
% \, %\forall \a \in H_1(\Sigma,\Z) \right\rbrace,$$ 
$$
\bar H(G)=\bigg\{  \omega \in H^{0,1}(k_\Sigma+D(G)) \, |\, Res_{r_i} \omega=0, 
 1 \leqslant i \leqslant \mu,
$$
$$\left. \,Re \int_\a
 (1-\g^2,i(1+\g^2),2\g)\omega=0, \, \forall \a \in H_1(\Sigma,\Z) \right\rbrace,$$

\noindent
where $k_\Sigma$ is a canonical divisor of $\Sigma$ and $H_1(\Sigma,\Z)$
is the first group of homology of $\Sigma.$
Suppose that the divisor $D$ %$2k_\Sigma+R(G)$ 
has an expression
of the form
$\sum n_j v_j -\sum m_i u_i,$ with $n_j,m_i \in \N.$
An element of  $H^{0,1}(D%2k_\Sigma+R(G)
)$ can 
be expressed as $f dz,$ where $f$ is
 a meromorphic function on $\Sigma$
with poles of order not bigger than $n_j$ at $v_i$ 
and zeroes of order not smaller than $m_i$ at $u_i.$ 
Equivalently, if 
 $g dz,$ where $g$ is a meromorphic function, 
is the  differential form associated with the 
divisor $D
%2k_\Sigma+R(G)
,$
the product $fg$ must be holomorphic.\\

\noindent
For $\omega \in \bar H(G),$ let
$X(\omega): \Sigma \setminus \{r_1,\ldots, r_\mu \} \to {\mathbb R}^3$ be the 
conformal immersion defined by
$$ X(\omega)(p)=Re \int^p  (1-\g^2,i(1+\g^2),2\g) \omega.$$%\frac{\sigma}{d\g}.$$
Then $X(\omega)\cdot G,$ the support function of $X(\omega),$ extends over 
the ramification points $r_1,\ldots,r_\mu$ smoothly and thus gives an element of $N(G).$
Conversely, every element of $N(G)$ is obtained in this way. In fact the map
\begin{equation}
\label{isomorfismo.teorema4}
\begin{array}{ccc}
i :\bar H(G) &\to& N(G)/L(G)\\
\omega & \to & [X(\omega)\cdot G ]\\
\end{array}
\end{equation}
is an isomorphism.
This result, used in association with the Weierstrass
 representation formula,
gives a  description of the space $N(G).$ 
To obtain the dimension of $N(G)$ it is sufficient to
 compute the dimension of $\bar H(G).$ Since the dimension 
of $L(G)$ is equal to 3, then ${\rm Nul}(G)=3+{\rm dim}\,\bar H(G).$\\

\noindent
We denote by $A_t$ a one parameter family ($0 <t<+\infty$)  of 
conformal diffeomorphisms
of the sphere $S^2$ defined by
$$\pi \circ A_t \circ \pi^{-1} w=tw,\quad w \in \C \cup \{\infty\}.$$
 
\noindent
We define for $0<t< \infty$
\begin{equation}
\label{Gt}
G_t=A_t \circ G.  
\end{equation}

\noindent
S. Nayatani in \cite{N2} gave lower and upper bounds 
for the index and, 
applying the method recalled above, for the nullity  of 
$G_t,$ $t\in (0,\infty),$ a deformation of an arbitrary 
holomorphic map  $G: \Sigma \to S^2,$
where $\Sigma$ is a compact Riemann surface.  In the same work,  
choosing appropriately the map $ G$ and the surface $\Sigma,$ 
he computed the index and 
the nullity for the Gauss map of  the Costa surface. In fact 
the extended Gauss map of this surface is a deformation of $G$ for 
a particular value of $t.$ We describe briefly the 
principal steps to get this result. \\

\noindent
Firstly it is necessary to study the vector space $\bar H( G_t).$ 
A differential  $\omega \in  H^{0,1}(k_\Sigma+D(G))$ with null residue
at the ramification points, is an element 
of $\bar H( G_t)$ if and only if the pair $(t \g,\omega)$ defines 
a branched minimal  surface by the Weierstrass representation. If 
one sets $\g=1/\wp'$  then there exist only 
two values of $t,$ denoted by $t'<t'',$ for which the 
condition above is verified and moreover
${\rm  dim}\, H(G_t)=1.$ In other words,
thanks to the characterization of the non linear elements of $N(G_t)$ by 
the isomorphism described by (\ref{isomorfismo}), if $t=t',t'',$  
${\rm Nul}(G_t)=4.$ As for the index, if $t=t',t''$ then
${\rm Ind}(G_t)=5.$ Since $G_{t''}$ is the 
extended Gauss map of 
the Costa surface, one can state:
\begin{theorem}
\label{nullindex.Costa}
Let $\bar G$ be the extended Gauss map of the 
Costa surface. Then
$${\rm Nul}(\bar G)=4, \quad  {\rm Ind}(\bar G)=5.$$
\end{theorem} 
 
\noindent
The same author in \cite{N1} treated the more difficult case of the 
Costa-Hoffman-Meeks surfaces of genus $2\leqslant g \leqslant 37$
by a slightly different method. 
That is the subject of next section.\\

%\noindent
%In the following section we shall restrict our attention to the Costa-Hoffman-Meeks surfaces
% and we will give some details of the computation in this case. We observe that the 
%extended Gauss map of this surface is $G_{t_2},$ where $G$ is a specific  
%holomorphic map and $t_2$ is one the values for which ${\rm Nul}(G_t)=4$ like 
%we will explain later. \\

\section{The case of the Costa-Hoffman-Meeks surface of genus smaller than $38$}
\label{caso.minore38}
In this section we expose some of the background details 
at the base of  section 3
of the work \cite{N1}. S. Nayatani provided them to us in \cite{N3}.\\

\noindent
We denote by $M_g$ the Costa-Hoffman-Meeks surface of genus $g.$
Let $\Sigma_g$ be the compact Riemann surface 

\begin{equation}
\label{Sigma.g}
\Sigma_g=\{ (z,w)\in (\C \cup \{\infty\})^2 \,|\, w^{g+1}=z^g (z^2-1)  \}
\end{equation}

\noindent
and let $Q_0=(0,0),P_+=(1,0),P_-=(-1,0),P_\infty=(\infty, \infty).$ It is 
known that $M_g=\Sigma_g \setminus \{P_+,P_-,P_\infty\}.$\\

\noindent
The following result describes the properties of symmetry
of $M_g$ and $\Sigma_g.$
\begin{lemma}(\cite{HM0})
\label{lemma.simmetrie}
Consider the conformal mappings of $(\C \cup \{ \infty\})^2:$
\begin{equation}
\label{lambdakappa}
\kappa(z,w)=(\bar z, \bar w) \quad \lambda(z,w)=(-z,\rho w),
\end{equation}
where $\rho=e^\frac{i \pi g}{g+1}.$  The map
$\kappa$ is of order $2$ and $\lambda$ is of order $2g+2.$ The group
generated by $\kappa$ and $\lambda$ is the dihedral group 
$D_{2g+2}.$ This group of conformal diffeomorphisms leaves $M_g$ 
invariant, fixes both $Q_0$ and $P_\infty$ and extend to $\Sigma_g.$
Also $\kappa$ fixes  
the points $P_\pm$ while $\lambda$ interchanges them.
\end{lemma}

%\begin{corollary}
%\label{corollario.simmetrie}
%The  surface $M_g$ possesses a dihedral group $D_{2g+2}$
%of conformal diffeomorphisms that fix $Q_0$ and extend 
%to $\Sigma_g.$ An element of $D_{2g+2}$ is of the form 
%$\tau=\lambda^j \kappa^i, $ for $i=0,1,$ $j=0,\ldots,2g+1.$
%If $j$ is even, $\tau$ fixes the points $P_\pm;$ if
%$j$ is odd, $\tau$ interchanges them. 
%\end{corollary} 

\noindent
We set  $\gamma(w)=w.$ Let $G:\Sigma_g \to S^2$ be the  holomorphic map defined by 
\begin{equation}
\label{definizione.G}
\pi \circ G(z,w)=\g(w).
\end{equation}

\noindent
%Let $\g$ be the meromorphic function defined by $\pi \circ G.$ 
%Let $p_j$ be the poles of $g$ and $r_i$ its ramification points.
We  denote by $r_i,$ $i=1,\ldots,\mu,$ the ramification points of $\g$
and by $R(G)$ the ramification divisor $\sum_{i=1}^{\mu} r_i.$ 
%(see subsection \ref{RiemannRoch}).
%, using theorem 5 of \cite{MR},
Theorem 5 of \cite{MR} shows that the space $N(G)/L(G),$ 
that we have introduced in previous section, is  also 
%be described in terms of meromorphic quadratic differentials.
%in terms of differential $1$-forms given in previous section,
%now, following theorem 5 of \cite{MR}, we will show 
%that the same space is 
isomorphic to a space of meromorphic quadratic differentials.
This alternative description  of $N(G)/L(G)$ that we present 
in the following, was adopted by S. Nayatani in \cite{N1}. 
We start defining the vector spaces $\hat H(G)$ and $H(G).$

\begin{equation}
\label{Hcappello.G}
\hat H(G)= \left\lbrace  \sigma \in H^{0,2}(2k_\Sigma+R(G)) \,| \, Res_{r_i} \frac{\sigma}{d\g}=0, i=1,\ldots,\mu
\right\rbrace,
\end{equation}

$$H(G)=\left\lbrace  \sigma \in \hat H(G) \, | \,Re \int_\a (1-\g^2,i(1+\g^2),2\g)\frac{\sigma}{d\g}=0, \, \forall \a \in H_1(\Sigma,\Z) \right\rbrace,$$ 

\noindent
where $k_\Sigma$ is a canonical divisor of $\Sigma.$ 
We remark that the elements of $H^{0,2}(2k_\Sigma+R(G))$ are quadratic differentials
(see subsection \ref{RiemannRoch}).
Since hereafter we will work only with quadratic differentials,
we can set $H^0(\cdot)=H^{0,2}(\cdot)$ to simplify the notation. 
If we suppose that the divisor $2k_\Sigma+R(G)$ 
has an expression
of the form
$\sum n_j v_j -\sum m_i u_i,$ with $n_j,m_i \in \N,$
an element of  $H^0(2k_\Sigma+R(G))$ can 
be expressed as $f (dz)^2,$ where $f$ is
 a meromorphic function on $\Sigma$
with poles of order not bigger than $n_j$ at $v_i$ 
and zeroes of order not smaller than $m_i$ at $u_i.$ 
Equivalently, if 
 $g (dz)^2,$ where $g$ is a meromorphic function, 
is the  differential form associated with the 
divisor $2k_\Sigma+R(G),$
the product $fg$ must be holomorphic.\\

\noindent
For $\sigma \in  H(G),$ let
$X(\sigma): \Sigma \setminus \{r_1,\ldots, r_\mu \} \to {\mathbb R}^3$ be the 
conformal immersion defined by
$$ X(\sigma)(p)=Re \int^p  (1-\g^2,i(1+\g^2),2\g) \frac{\sigma}{d\g}.$$
Then $X(\sigma)\cdot G,$ the support function of $X(\sigma),$ extends over 
the ramification points $r_1,\ldots,r_\mu$ smoothly and thus gives an element of $N(G).$
Conversely, every element of $N(G)$ is obtained in this way. In fact the map
\begin{equation}
\label{isomorfismo}
\begin{array}{ccc}
i :H(G) &\to& N(G)/L(G)\\
\sigma & \to & [X(\sigma)\cdot G ]\\
\end{array}
\end{equation}
is an isomorphism.
%This result, used in association with the Weierstrass
% representation formula,
%gives a holomorphic description of the space $N(G).$ 
So to obtain the dimension of $N(G)$ it is sufficient to
 compute the dimension of $H(G).$ We recall that the dimension 
of $L(G)$ is equal to 3, so ${\rm Nul}(G)=3+{\rm dim}\,H(G).$\\

\noindent
Since the extended Gauss map of the Costa-Hoffman-Meeks
surfaces is a deformation in the sense of the definition (\ref{Gt})
of the map $G,$ we need to study the space $H(G_t).$ 
From  (\ref{Hcappello.G}) and (\ref{Gt}) it is clear 
that $\hat H(G)=\hat H(G_t)$ and  
% $$H(G_t)=\left\lbrace  \omega \in \hat H(G_t) \, | \,Re \int_\a (1-t^2\g^2,i(1+t^2 \g^2),2t\g)\o=0, \, \forall \a \in H_1(\Sigma,\Z) \right\rbrace.$$ 
$$H(G_t)=\left\lbrace  \sigma \in \hat H(G_t) \, | \,Re \int_\a (1-t^2\g^2,i(1+t^2 \g^2),2t\g)\frac{\sigma}{d\g}=0, \, \forall \a \in H_1(\Sigma_g,\Z) \right\rbrace.$$ 

\noindent
Long computations (\cite{N3}, see subsection \ref{base} for some details) show 
that a basis of the differentials of the form  $\sigma/d\g,$
where $\sigma \in \hat H(G)=\hat H(G_t),$ and
whose residue at the ramification points of $\g(w)=w$ is zero,
is formed by 
$$ \omega_k^{(1)}= \frac {z^{k-1}}{w^k}  \frac{dz}{w},
\quad {\rm with}\quad k=1,\ldots, g-1,$$
$$\omega_k^{(2)}= \frac{((k-2) z^2 -k A^2)}{ (z^2-A^2)^2}\left(\frac {z}{w}\right)^{k-1} \frac{dz}{w},\quad{\rm with}\quad k=0,\ldots, g,$$ 
$$ \omega_k^{(3)}= \frac{((k-2) z^2 -k A^2)}{w (z^2-A^2)^2}\left(\frac {z}{w}\right)^{k-1} \frac{dz}{w},\quad {\rm with}\quad k=0,\ldots, g-1,$$ 
where $A = \sqrt{\frac g {g+2}}.$\\

\noindent
Now we put attention to the space $H(G_t).$
We recall that we are interested in the computation
of its dimension.
By the definition of $H(G_t),$ a differential 
 $\sigma \in \hat H(G_t)$ belongs to $H(G_t)$ if 
and only if $\forall \a \in H_1(\Sigma_g,\Z)$ the 
differential form 
$\omega=\frac{\sigma}{d\g}=\frac{\sigma}{dw}$
satisfies
\begin{equation}
\label{condizione1}
\int_\a \omega= t^2  \overline{ \int_\a \g^2(w) \omega },
\end{equation}
\begin{equation}
\label{condizione2}
 Re \int_\a \g(w)  \omega =0.
\end{equation}
If these two conditions are satisfied then
$(\g,w)$ are the Weierstrass data of a branched minimal surface. 
Of course, it is sufficient to impose that these equations are satisfied 
when $\a$ varies between the elements of a basis 
of $H_1(\Sigma_g,\Z).$ The convenient basis of 
$H_1(\Sigma_g,\Z)$ is constructed as follows.
Let $\b(s)=\frac 1 2 +e^{i 2 \pi s},$ 
$0\leqslant s \leqslant 1. $
Let $\tilde \b(s)=(\b(s), w(\b(s)))$ be a lift of $\b$ 
to $\Sigma_g$ such that, for example, $\tilde \b(0)=(\frac 3 2,w(0)), $
 with $w(0) \in \R.$
As stated in lemma \ref{lemma.simmetrie} 
the group  of conformal diffeomorphisms of $\Sigma_g$ is
isomorphic to the dihedral group $D_{2g+2}.$
The collection $\{\lambda^l \circ \tilde{\b}, l=0,\ldots,2g-1\},$
where $\lambda$ is the generator of $D_{2g+2}$
of order $2g+2,$ is a basis of $H_1(\Sigma_g,\Z)$
(see \cite{HM0}).\\

\noindent
Now we must impose
 (\ref{condizione1}) and (\ref{condizione2}) 
for $\a=\lambda^l\circ \tilde \b,$ with $l=0,\ldots, 2g-1.$ 
 To do that we collapse $\b$ to the unit interval.
 In other terms we deform continously $\b$ in such a way
 the limit curve is the union of two line segments 
 lying on the real line.    
 We set  $$\omega =\sum_0^{g-1} c_k^{(1)} \omega_k^{(1)} +
\sum_0^{g} c_k^{(2)} \omega_k^{(2)}
+\sum_0^{g-1} c_k^{(3)} \omega_k^{(3)},$$
where $c_k^{(i)} \in \C.$\\

\noindent
Taking into account  these assumptions, it is possible
to show that the equation (\ref{condizione1}), if the genus $g$ is 2, 
is equivalent  to the following  system of four equations
(see subsection \ref{sistemi})
 \begin{equation}
\label{sistema1}
\left\{ \begin{array}{ll}
f_0=-t^2 \bar h_{0}\\
f_1=0\\
p_{1}=-t^2 \bar q_{1}\\
p_{2}=-t^2 \bar q_{0}.\\
\end{array} \right.
\end{equation}

\noindent
If $g\geqslant 3$ 
there are the following  additional $2g-4$ equations to consider
 \begin{equation}
\label{sistema1.++}
\left\{ \begin{array}{ll}
%f_0=-t^2 \bar h_{0}\\
%f_1=0\\
f_{k}=-t^2 \bar q_{g-k+2}\\
%p_{1}=-t^2 \bar q_{1}\\
%p_{2}=-t^2 \bar q_{0}\\
p_{g-k+2}=-t^2 \bar h_{k}\\
\end{array} \right.
\end{equation}
where $k=2,\ldots,g-1$ and

$$ f_0= \frac{(g+2)^2}{2(g+1)} \,c_0^{(3)} \sin \left( \frac{\pi}{g+1}\right)  K_0,$$

$$ f_k= \left( -c_k^{(1)} +\frac{(g+2)(g+2+k)}{2(g+1)} c_k^{(3)} \right)
 \sin \left( \frac{(k+1)\pi}{g+1}\right)  K_k, \quad k=1,\ldots,g-1,$$

$$ h_0= \frac{(g+2)^2}{2(g+1)}\; c_0^{(3)} \sin \left( \frac{-\pi}{g+1}\right)  J_0,$$

$$ h_k= \left( c_k^{(1)} +\frac{(g+2)(g+2-k)}{2(g+1)}\, c_k^{(3)} \right)
 \sin\left(  \frac{(k-1)\pi}{g+1}\right)  J_k, \quad k=2,\ldots,g-1,$$

$$ p_k= -\frac{(g+2)k}{2(g+1)}\, c_k^{(2)} \sin \left( \frac{k \pi}{g+1}\right)  I_k,
 \quad k=1,\ldots,g,$$

$$q_k=\frac{(g+2)(2g+4-k)}{2(g+1)} \,c_k^{(2)}\sin \left( \frac{(k-2)\pi}{g+1}\right)  L_k,
 \quad k=0,1,3,\ldots,g,$$
and

$$ I_m= \frac{g+1}{m} \frac{\G\left( 1+ \frac{m}{2(g+1)}\right) \G\left( 1- \frac{m}{g+1} \right) } {\G\left( 1- \frac{m}{2(g+1) }\right)}, $$

$$ J_m= \frac{g+1}{g-m+2} \frac{\G\left( \frac 12+ \frac{m-1}{2(g+1)}\right) \G\left( 1- \frac{m-1}{g+1} \right) } {\G\left( \frac12- \frac{m-1}{2(g+1) }\right)}, $$

$$ K_m= J_{m+2},$$

$$ L_m= \frac{m-2}{2g-m+4} I_{m-2}.$$
The equation (\ref{condizione2}) if the genus $g$ is 2, 
is equivalent to the following system of two equations
(see subsection \ref{sistemi})
 \begin{equation}
\label{sistema2}
\left\{ \begin{array}{ll}
d_1=0\\
e_{2}=\bar e_{0}.\\
\end{array} \right.
\end{equation}
If $g\geqslant  3$ there are the following  additional $g-2$ equations to consider
 \begin{equation}
\label{sistema2.++}
%\left\{ \begin{array}{ll}
%d_1=0\\
d_{k}= \bar e_{g-k+2}\\
%e_{2}=\bar e_{0}\\
%\end{array} \right.
\end{equation}
where $k=2,\ldots,g-1,$ and

$$d_k=\left( c_k^{(1)}-\frac{k(g+2)}{2(g+1)}c_k^{(3)}\right) 
\sin \left( \frac{k\pi}{g+1}\right)  I_k, \quad k=1,\ldots,g-1,$$

$$e_k= \frac{(g+2)(g+2-k)}{2(g+1)} c_k^{(2) }
\sin \left( \frac{(k-1)\pi}{g+1}\right)  J_k, \quad k=0,2,\ldots,g.$$

%\noindent
%We observe that if $g=2$ this system reduces to 
% $$
%\left\{ \begin{array}{ll}
%f_0=-t^2 \bar h_{0}\\
%f_1=0\\
%p_{1}=-t^2 \bar q_{1}\\
%p_{2}=-t^2 \bar q_{0}\\
%\end{array} \right.
%$$
\noindent
We are looking for the values of $t$ such that the previous systems 
  have non trivial solutions in terms of  $c_i^{(j)}.$ Only for
   these special values of $t$ it  holds  $\dim H(G_t)> 0$ or
equivalently  ${\rm Nul}(G_t)>3.$\\

\noindent 
We start with the analysis of the system (\ref{sistema1}). 
This system admits  non trivial solutions  if and only if $t$ takes 
 three values denoted by  $t_1,t_2,t_3.$
Obviously they are functions of $g.$\\

\noindent
If we set $s= \frac 1 {g+1}$ then we can write
$$t_1=\sqrt{\frac{K_0}{J_0}}=  
%\sqrt{\frac{g(g+2)}{4(g+1)^2}}
\frac{\sqrt{1-s^2}}{2}
\sqrt{\frac{\G\left(  1-s \right) }{\G\left(  1+s \right)}}
\frac{\G\left(  1-\frac s {2} \right) }{\G\left(  1+\frac s {2} \right)},$$
$$t_2=\sqrt{\frac{I_1}{(2g+3)L_1}}=
\sqrt{\frac{\G\left(  1- s \right) }{\G\left(  1+  s \right)}}
\frac{\G\left(  1+\frac s {2} \right) }{\G\left(  1-\frac s {2} \right)},$$
$$t_3=\sqrt{\frac{I_2J_0}{gL_0K_0}}=
\frac{2}{1-s} \sqrt{\left( \frac{\G(1+s)}{\G(1-s)}\right)^ {3}}
 \sqrt{ \frac{\G(1-2s)}{\G(1+2s)}}\,\frac{\G(3/2-s/2)}{\G(1/2+s/2)}.$$

\noindent
We recall that if $g\geqslant 3$ there are other equations to consider.
They are 
$$
\left\{ \begin{array}{ll}
f_{k}=-t^2 \bar q_{g-k+2}\\
p_{g-k+2}=-t^2 \bar h_{k}\\
d_{k}= \bar e_{g-k+2}\\
\end{array} \right.
$$
where $k=2,\ldots,g-1.$ 
Thanks to the particular structure of the equations, it is possible
to study separately for each set of three equations 
the existence of solutions. 
Each set of three equations  admits non 
trivial solutions if and only if the following matrix has 
determinant equal to zero
 $$
\left( \begin{array}{ccc}
-K_k &   (g+2+k)  K_k&   (g+2+k)t^2 L_{g-k+2} \\
 t^2 J_k&   (g+2-k)  t^2J_k&  (g+2-k) I_{g-k+2}  \\
I_k & -k I_k& -k J_{g-k+2}\\
\end{array} \right ).
$$
After  the change of variable $l=g-k+1$ so that $2 \leqslant  l \leqslant g-1,$ 
it is possible to show that the determinant is 
\begin{equation}
-(g+2)(at^4 +b t^2+ c),
\end{equation}
with
$$a= (2g-l+3) I_{g-l+1}J_{g-l+1} L_{l+1}$$
$$b= -2 (g-l+1)  J_{l+1}  J_{g-l+1} K_{g-l+1} $$
$$c=(l+1) I_{g-l+1} I_{l+1}K_{g-l+1} .$$\\
%where $2 \leqslant l \leqslant g-1.$\\
\noindent
We are interested in finding the positive values of $t$ such that 
\begin{equation}
\label{equationt}
at^4 +b t^2+ c=0.
\end{equation}

\noindent
%In fact he has showed that the discriminant of \ref{equationt} is an increasing
%function. 
To simplify the notation we introduce the following three functions
$$F(v)=
\left(\frac{\G(\frac{1}{2}+\frac{v}{2})}{\G(\frac{1}{2}-\frac{v}{2})}\right)^ {2}
\frac{\G(1-v)}{\G(1+v)},$$

$$I(v)=
\left(\frac{\G(1-\frac{v}{2})}{\G(1+\frac{v}{2})}\right)^ {2}
\frac{\G(1+v)}{\G(1-v)},$$

$$L(v)=\left(\frac{\G(1+\frac{v}{2})}{\G(1-\frac{v}{2})}\right)^ {2}
\frac{\G(1-v)}{\G(1+v)}=\frac{1}{I(v)}.  $$

\noindent
The discriminant $b^2-4ac$ of the equation (\ref{equationt}), seen like an equation of 
degree two in the variable $t^2,$ is negative 
if and only if $X=b^2/4ac <1.$ 
It is possible to show that   
\begin{equation}
\label{determinante}
X= \frac{l^2} {l^2-1} F^2\left (\frac l {g+1} \right)  
I\left (\frac {l-1} {g+1} \right) I\left (\frac {l+1} {g+1} \right).
\end{equation}
 
\noindent
S. Nayatani %in \cite{N3} 
showed that if $2\leqslant g \leqslant 37,$ then 
$X<1$ and as consequence the equation (\ref{equationt}) 
has not any solution   since  its discriminant is negative.  
Then %if $2\leqslant g \leqslant 37$ 
$\dim H(G_t)>0$ only for 
$t=t_1,t_2,t_3.$ Summarizing we can state
(see \cite{N1} for other details):

\begin{theorem}
\label{nullity.37}
If $2\leqslant g \leqslant 37$ and $t \in (0,+\infty),$ then
%As consequence ${\rm Nul}(G_t)=3$
%for $t \in (0, +\infty)$ excepted when $t=t_1,t_2,t_3$ 
%for which we have 

$$
{\rm Nul}(G_t)=\left\lbrace  
\begin{array}{cl}
4 & \hbox{if}\quad t=t_1,t_2\\
5 & \hbox{if}\quad t=t_3\\
3 & \hbox{elsewhere.} \\
\end{array}
\right.$$
\end{theorem}
\noindent
Since the extended Gauss map of  the Costa-Hoffman-Meeks 
surfaces is exactly $G_{t_2},$
it is possible to state  that the null space of the Jacobi operator
of $M_g$ has dimension equal to $4$ for $2 \leqslant g \leqslant 37.$\\

\noindent
Other values of $t$ for which ${\rm Nul}(G_t)>3$ are admitted only if $g\geqslant 38.$
In \cite{N1} S. Nayatani conjectured these values were bigger than $t_3.$ 
The proof of the conjecture and its consequences will be showed in  sections \ref{caso.38} and \ref{statements}.\\

\section{The case $g \geqslant 38$}
\label{caso.38}

\noindent
%Thanks to the previous observations it is clear
%that if we assume $g \geqslant 38$ 
S.Nayatani proved 
%It is possible to prove %(\cite{N3}) 
that $X$ 
%(and so the discriminant of \ref{equationt}) 
is a decreasing function in the variables $l,$
$$x=\frac{l}{g+1}, y= \frac{l+1}{g+1}, z= \frac{l-1}{g+1}$$
with $2 \leqslant l \leqslant g-1.$
We recall that we  have set $s=\frac{1}{g+1}.$
We know that for $l=2$ and $g=37$ the discriminant of the 
equation (\ref{equationt}) is negative.
For these values of $l$ and $g$ the variables $x,y,z,s$ are 
respectively equal to $x_{max}=2s_{max},$
$y_{max}=3s_{max},$ $z_{max}=s_{max}=1/38.$
Then we will study the solutions of (\ref{equationt}) for $i \in [0,i_{max}]$
(we call admissible values the values in these intervals ) 
where $i$ denotes $x,y,z,s,$ 
because for bigger values of the  variables the discriminant continues to 
be negative and so the 
equation (\ref{equationt}) does not admit solutions. \\
%Numerical tests show that %the maxima of the admissible values,
%the value $i_{max},$  become smaller as $g$ is bigger. Since it is not possible
%to explicit the dependence of $i_{max}$
%on $g$ we shall work with constant quantities.\\

\noindent
All the solutions of (\ref{equationt}), that we denote by $t_\pm(l,g),$ satisfy
%$$F(x)= {\left( \frac{\G(\frac{1}{2}+\frac{s}{2})}{\G(\frac{1}{2}-\frac{s}{2})}\right)^ {2}
 %\frac{\G(1-x)}{\G(1+x)}}.$$
%$$G(s)= \left(  \frac{\G(1-\frac s 2)}{\G(1 + \frac s 2)}\right)^ {2}
% {\frac{\G(1+s)}{\G(1-s)}}$$
%$$P(l,g)= \frac{l}{l-1} F \left( \frac{l}{g+1}\right)  G\left( \frac{l-1}{g+1}\right) $$
%$$Q(l,g)= \frac{l+1}{l-1} L \left( \frac{l+1}{g+1} \right)  G\left( \frac{l-1}{g+1}\right)$$
$t_\pm^2(l,g)= T_1\pm T_2,$ with
\begin{equation}
\label{primo.addendo}
%T_1=\frac{l}{l-1} F \left( \frac{l}{g+1}\right)  I\left( \frac{l-1}{g+1}\right)
T_1=\frac{l}{l-1} F \left( x\right)  I\left( z\right)
\end{equation}
 and 
\begin{equation}
\label{secondo.addendo}
%T_2=\sqrt{\left( \frac{l}{l-1}\right)^2 F^2 \left( \frac{l}{g+1}\right)  I^2\left( %\frac{l-1}{g+1}\right)- \frac{l+1}{l-1} L \left( \frac{l+1}{g+1} \right)  I\left( %\frac{l-1}{g+1}\right)}.
T_2=\sqrt{\left( \frac{l}{l-1}\right)^2 F^2 \left( x\right)  I^2\left(z\right)- \frac{l+1}{l-1} L \left(y \right)  I\left( z\right)}.
\end{equation}

\noindent
We will prove that, for all the values of $l$ and $g,$ such that $0\leqslant \frac{l}{g+1}\leqslant x_{max}=\frac{2}{38},$ 
with  $2 \leqslant l \leqslant g-1$ and $g \geqslant 38,$ such that 
$T_2$ is a real number, it holds

\begin{equation}
\label{inequality}
t_3^2 \left(s \right)< t_-^2(l,g).
\end{equation}

\noindent
We need  study the behaviour of the functions $F,I,L,F^2,I^2$ that 
appear in (\ref{primo.addendo}) and (\ref{secondo.addendo}).
This aim is pursued by the use of zero order series of these
functions.

\noindent
The Mac-Laurin  series  of  the functions $F(x),G(z),L(y),F^2(x),I^2(z)$ for 
admissible values of $x,y,z$ are 
%$$F(x)= 1+ 2\left(\g+\psi\left( \frac{1}{2}\right)\right )x+
%2\left(\g+\psi\left( \frac{1}{2}\right)\right )^2x^2
%+R(x) x^3$$
\begin{equation}
\label{sviluppo.F}
F(x)= 1+ R_F(d_1 x) x, \quad I(z)=1+R_I(d_2 z) z ,\quad L(y)=1+R_L(d_3 y) y ,
\end{equation}
 %$$$$
 %\quad -0.055 \leqslant R_I(z) \leqslant 0,$$
%where $-0.055 \leqslant R_I(z) \leqslant 0,$\\
%$$$$
%\quad 0 \leqslant R_L(y) \leqslant 0.055.$$
%where $0 \leqslant R_L(y) \leqslant 0.055,$ 
$$F^2(x)= 1+R_{F^2}(c_1 x) x ,\quad I^2(x)= 1+R_{I^2}(c_2 x) x ,$$
%\quad 13\leqslant R_{F^2}(x) \leqslant 15$$
%where $ 13\leqslant R_{F^2}(x) \leqslant 15.$
%$$$$
%\quad -0.11 \leqslant R_{I^2}(x) \leqslant 0$$
%where $-0.11 \leqslant R_{I^2}(x) \leqslant 0,$
where $c_i,d_i\in(0,1).$ 
So we can write
$$F (x)  I (z)=  1+ R_{FI}(x,z), \,F^2 (x)  I^2 (z)=  1+ R_{F^2I^2}(x,z),\,
L (y)  I (z)=  1+ R_{LI}(y,z),$$ 
\noindent
with $$R_{FI}(x,z)=R_F(d_1 x) x +R_I(d_2 z) z+ R_F(d_1 x) R_I(d_2 z) x z,$$
 $$R_{F^2I^2}(x,z)=R_{F^2}(c_1 x) x +R_{I^2}(c_2 z) z+ 
R_{F^2}(c_1 x) R_{I^2}(c_2 z) x z,$$
 $$R_{LI}(y,z) =R_L(d_3 y) y+R_I(d_2 z) z+R_I(d_2 z) R_L(d_3 y) z y.$$\\
%\noindent
%Using these series we start giving a rough estimate of the quantity $ l^2 x=
%\frac {l^3} %{g+1}.$
%\begin{lemma}
%\label{stima.l2x}
%$$ l^2 x=\frac {l^3} {g+1}< %\frac{1}{8 \ln2}
%1 .$$
%\end{lemma}
%\noindent
%{\bf Proof.} 
%We consider the expression of the function $Q$ given by (\ref{determinante}).
%Let us impose the condition $Q>1$ using an approximate expression for $Q.$
%We find $$Q> \frac{l^2}{l^2-1}(1+2Ax)>1.$$
%From this it follows that 
%$$ l^2 x=\frac {l^3} {g+1}<-\frac 1 {2A}=\frac{1}{8 \ln2} \sim 0.180.$$
%\hfill \qed
%Because of the use of an approximated  expression we don't obtain the correct value
%($\sim 0.191$ using numerical tests), but it's sufficient for our aim. 
%Moreover taking into account the relations between $x,y,z,$ it is easy to get 
%similar estimates for $l^2y=l^2x+l/(g+1)$ and $l^2z=l^2x-l^2/(g+1).$ \\
\noindent
In the following $\psi(x)$ denotes  the digamma function. It is related 
to $\G(x),$ the gamma function,  by  
$$ \psi(x)= \frac{d}{dx} \left( \ln \G(x)  \right). $$
For the properties of these special functions we will refer to \cite{AS}.\\
%\noindent
%%and $\gamma=-\psi(1)$ the constant of Euler-Mascheroni.
%We denote $C=2\left(-\psi(1)+\psi\left( \frac{1}{2}\right)\right)=-4 \ln 2.$\\

\noindent
The following proposition gives useful properties of the functions 
just introduced.
\begin{proposition}
\label{proprietà}
If $x \in [0,x_{max}], $ $z \in [0,z_{max}]$ and  $y \in [0,y_{max}],  $ 
the following assertions hold:
\begin{enumerate}
\item \label{proprietà.RF} $R_{F}(x)<0$ %is a negative function.
\item \label{proprietà.RI} $R_{I}(z)\leqslant 0$ %is not positive and decreasing. 
\item \label{proprietà.RL}  $R_{L}(y)\geqslant0$ %is not negative and increasing.
\item \label{proprietà.RF'} $(R_{F})_x'(x)>0$
\item \label{proprietà.RI'} $min (R_{I})_z'(z)= -0.095\cdots$
\item \label{proprietà.RFI}   $R_{FI}(x,z)\geqslant Cx$ with $C=-4 \ln 2$
\item \label{proprietà.RLI}   $R_{LI}(y,z)\geqslant 0$% is not negative.
\item \label{proprietà.RI2}  $R_{I^2}(z)\leqslant 0$% is not positive and decreasing.
\item \label{proprietà.RF2} $W(x)=R_{F^2}(x)<0$ %is negative.
\item \label{proprietà.derivata.1W} $W'_x(x) >0$, so $R_{F^2}(x)$ is 
an increasing function
\item \label{proprietà.derivata.2W}   $W''_{xx}(x)<0$ %is negative.
\item  \label{proprietà.derivata.3W}  $W'''_{xxx}(x)>0$% is positive.
\item  \label{proprietà.derivata.1Y} If we set $Y(x)=x\, W(x),$ then 
$Y'_{x}(x)< 0$ 
\item  \label{proprietà.derivata.2Y} $Y''_{xx}(x)> 0$
\item  \label{proprietà.derivata.3Y} $Y'''_{xxx}(x)< 0.$

\end{enumerate}

\end{proposition}
\noindent
{\bf Proof.} 
\begin{enumerate}
\item $R_{F}(x)=F'_x(x)=F(x) \Psi_F(x),$ where
$$  \Psi_F(x)=-\psi(1-x)-\psi(1+x)+\psi\left(\frac 1 2 -\frac x 2 \right)+
\psi\left (\frac 1 2+\frac x 2 \right).  $$
We observe that $$\Psi_F(x)= %2\psi\left(\frac 1 2 \right)- 2\psi(1)+
2 \sum_{k=0}^\infty \frac{1}{(2k)!} \left(\frac{1}{2^{2k}} \psi^{(2k)}\left(\frac 1 2 \right)  - \psi^{(2k)}(1)\right) x^{2k}.$$
Since  $\Psi_F(0)=2\psi\left(\frac 1 2 \right)- 2\psi(1)=-4\ln2,$  
$\psi^{(2k)}(1)<  0$  and  
$\psi^{(2k)}\left(\frac 1 2 \right)= (2^{2k+1}-1) \psi^{(2k)}(1)<0,$ 
if $k\geqslant 1$ (see formulas  6.4.2 and 6.4.4 of  \cite{AS}), 
we can conclude that $\Psi_F(x)<0$ and it  is a decreasing function. 
Since $F(x)>0$    then $R_{F}(x)<0$ and $F(x)$ is a decreasing function.\\

\item $R_{I}(z)=I'_z(z)=I(z) \Psi_I(z),$  where

$$  \Psi_I(z)=\psi(1-z)+\psi(1+z)-\psi\left(1-\frac z 2 \right)-
\psi\left (1+\frac z 2 \right).$$
We observe that $$\Psi_I(z)=2\sum_{k=1}^\infty
\frac{1}{(2k)!} \psi^{(2k)}(1)\left(1- \frac{1}{2^{2k}} \right) z^{2k}. $$
Since $\psi^{(2k)}(1)<0$ for $k\geqslant 1$  
then  $\Psi_I(z)\leqslant 0$ and it is  
a decreasing function. Since $I(z)>0$  then   $R_{I}(z)\leqslant 0.$ \\
%and it is a decreasing function. Since $G(0)=1,$ 
%We have also showed that $G(t)$ is  a positive and decreasing function.

\item $R_{L}(y)=L'_y(y)=L(y) \Psi_L(y),$  where $\Psi_L(y)=-\Psi_I(y).$
Then  $\Psi_L(y)\geqslant 0$ and it is an increasing function. 
Since $L(y)=1/I(y)>0,$  then   $R_{L}(y)\geqslant 0.$ \\
%and it is an increasing function. Since $L(0)=1,$ we have also showed 
%that $L(t)$ is a positive and increasing function.
\item The derivative of $R_F$ is $F''_{xx}(x)=F(x) (\Psi_F^2(x)+(\Psi_F)'_x(x)).$ 
Since $\Psi_F(x)<0$ and it is a decreasing function,  $\Psi_F^2(x)>0$ and  
increasing. It holds $\Psi_F^2(x)\geqslant \Psi_F^2(0)=16 \ln^2 2.$
$$(\Psi_F)'_x(x)=2 \sum_{k=1}^\infty \frac{1}{(2k-1)!} \left( \frac{1}{2^{2k}} \psi^{(2k)}\left(\frac 1 2 \right)  - \psi^{(2k)}(1)\right) x^{2k-1}. $$
All the coefficients of the series are negative (see the point 
\ref{proprietà.RF}) so $(\Psi_F)'_x(x)\leqslant 0$ and it is a decreasing 
function. In particular $(\Psi_F)'_x(x) \geqslant (\Psi_F)'_x(x_{max})= -0.19\cdots.$ 
Since $F(x)>0$ and it is a decreasing function we can conclude that
 $$F''_{xx}(x)\geqslant F(x_{max})(\Psi_F^2(0)+(\Psi_F)'_x(x_{max}) )=
  6.4\cdots.$$

\item The derivative of $R_I$ is  $I''_{zz}(z)=I(z) (\Psi_I^2(z)+(\Psi_I)'_z(z)).$ 
Since $\Psi_I(z)\leqslant 0$ and it is a decreasing function (see the point 
\ref{proprietà.RI}),  $\Psi_I^2(z)\geqslant 0$ and 
increasing. It holds 
$\Psi_I^2(z)\leqslant \Psi_I^2(z_{max})= 1.5\cdots \cdot 10^{-6}.$
$$(\Psi_I)'_z(z)=2\sum_{k=1}^\infty
\frac{1}{(2k-1)!} \psi^{(2k)}(1)\left(1- \frac{1}{2^{2k}} \right) z^{2k-1}. $$
All the coefficients of the series are negative 
 so $(\Psi_I)'_z(z)\leqslant 0$ and it is a decreasing 
function. In particular $(\Psi_I)'_z(z) \geqslant (\Psi_I)'_z(z_{max})= -0.095\cdots.$ 
Since $I(z)>0$ and it is a decreasing function we can conclude that 
 $$I''_{zz}\geqslant I(z_{max}) (\Psi_I^2(0)+(\Psi_I)'_z(z_{max}))
  =-0.095\cdots.$$

\item Since $R_F<0$ and $R_I\leqslant 0,$  it holds that
$$R_{FI}(x,z)\geqslant R_F(d_1 x) x +R_I(d_2 z) z,$$ where $d_i \in (0,1).$
The point \ref{proprietà.RF'} implies that 
$R_F$  is an increasing function and we have computed the positive minimum 
(that we denote by $m$) value of its derivative.
Thanks to the  point \ref{proprietà.RI'} we have
 $m >|n|,$ where $n$ denotes the negative minimum value
of the derivative of $R_I.$
%Now it is sufficient to remember that the variables $x$ and $z$ are
%not independent. 
%We can conclude that $R_{FI}$ is an increasing function. 
%NUOVO 13 febbraio
Now we observe that
$$R_F(d_1x)x+R_I(d_2z)z \geqslant (R_F(0)+m x)x +(R_I(0)+n z)z 
\geqslant$$
$$ R_F(0)x+R_I(0)z+(m+n)z^2 \geqslant R_F(0)x+R_I(0)z =Cx. $$ 
To obtain this chain of inequalities we used the fact that $m+n>0$
and $x\geqslant z.$ Then $R_{FI}\geqslant Cx.$\\

\item We recall that $R_{LI}(y,z)=L(y)I(z)-1,$ $L(t)=1/I(t)$ and  
$$y=\frac{l+1}{g+1} > \frac{l-1}{g+1} =z.$$
We want to prove that $L(y)I(z)-1\geqslant 0$ or equivalently
$L(y)\geqslant 1/I(z).$ But thanks to the point \ref{proprietà.RL},
 %$L$ is an increasing function, so 
we have
$$L(y)\geqslant L(z)=\frac{1}{I(z)}.$$ 

\item $R_{I^2}(z)=(I^2)'_z(z)=2 I^2(z) \Psi_I(z).$
From the proof of the point \ref{proprietà.RI}, $\Psi_I(z)\leqslant 0$ 
and it is a decreasing function. 
Since $2I^2(z)>0,$  then also $R_{I^2}(z)\leqslant 0.$\\
% and it is a decreasing function. We have also showed 
%that $I^2(z)$ is a decreasing function.

\item $W(x)=(F^2)'_x(x)=2 F^2(x) \Psi_F(x).$ 
In the point \ref{proprietà.RF} we have observed that $\Psi_F(x)$ is a negative and decreasing function. 
Since $2F^2(x)>0,$  then also $W(x)$ 
is a negative function.\\
% We have also showed that $F^2(t)$ is a positive and decreasing function.
\item %We set $W(t)=R_{F}(t).$ We want to prove that $W'_t(t)$ is positive.
$W'_x(x)= F^2\left( 4 \Psi_F^2(x)+2 (\Psi_F)'_x(x) \right).$
Since $\Psi_F(x)<0$ and it is a decreasing function,  $\Psi_F^2(x)$ is 
a positive and increasing function. In the proof of the 
point \ref{proprietà.RF'} we observed that 
%$$(\Psi_F)'_t(t)=2 \sum_{k=1}^\infty \frac{1}{(2k-1)!} \left( \frac{1}{2^{2k}} %\psi^{(2k)}\left(\frac 1 2 \right)  - \psi^{(2k)}(1)\right) t^{2k-1}. $$
%All the coefficients of the series are negative so 
$(\Psi_F)'_x(x)\leqslant 0$ and it is a decreasing 
function. %for $t\in I.$ 
Since $2 (\Psi_F)'_x(x_{max})= -0.38\cdots$ and 
$4 \Psi_F^2(x)\geqslant 4 \Psi_F^2(0) =64 \ln^2 2 = 30.74\cdots,$ we can conclude that 
$W'_x(x)>0.$\\
 
\item The explicit expression of $W''_{xx}$ is 

$$W''_{xx}=\frac{1}{2} F^2(x) 
\left( 16 \Psi_F^3(x) +24 \Psi_F(x)  (\Psi_F)'_x(x) +4(\Psi_F)''_{xx} (x) \right).  $$

\noindent
%We want to show that this function is negative.
In the proof of the point \ref{proprietà.RF} we  observed that $\Psi_F(x)$
is a negative and decreasing function. So $16 \Psi_F^3(x) \leqslant 16 \Psi_F^3(0)
% =16 C^3
=-1024 \ln^3 2= -341.\cdots.$ 
Thanks to the proof of the point \ref{proprietà.derivata.1W} we know that 
 $(\Psi_F)'_x (x)\leqslant 0$ and it is a decreasing function. In particular 
$0\geqslant (\Psi_F)'_x (x)\geqslant (\Psi_F)'_x (x_{max}) = -0.19\cdots.$ 
We can conclude that 
$$ 24 \Psi_F(x)  (\Psi_F)'_x (x) \leqslant 24  (\Psi_F)'_x (x_{max}) 
\Psi_F(x_{max})= 12.\cdots.$$ 
As for the last summand, it is negative. In fact 
$$(\Psi_F)''_{xx} (x)=  2 \sum_{k=1}^\infty \frac 1{(2k-2)!} 
\left( \frac{1}{2^{2k}} \psi^{(2k)}\left(\frac 1 2 \right)  -
\psi^{(2k)}(1)\right) x^{2k-2}.$$
Since all the coefficients of the series are negative, we get
$$4(\Psi_F)''_{xx} (x)\leqslant 4(\Psi_F)''_{xx} (0)=-12 \zeta(3)= -14.4\cdots,$$
where $\zeta(\cdot)$ denotes the Riemann zeta function.
\noindent
We can conclude that 
 $$16 \Psi_F^3(x) +24 \Psi_F(x)  (\Psi_F)'_x (x) +4  (\Psi_F)''_{xx} (x) \leqslant $$
%-341+12- 14.4=
$$\leqslant 16 \Psi_F^3(0) +24 \Psi_F(x_{max})  (\Psi_F)'_x (x_{max}) +
4  (\Psi_F)''_{xx} (0)=
-342.7\cdots .$$
That assures  $W''_{xx}<0.$\\

\item The explicit expression of $W'''_{xxx}$ is 

$$W'''_{xxx}=\frac{1}{4} F^2(x) 
\left( 64 \Psi_F^4 + 192 \Psi^2_F  (\Psi_F)'_x
 +48 ((\Psi_F)'_x)^2 +\right.$$
$$\left.+64  \Psi_F (\Psi_F)''_{xx}   +8  (\Psi_F)'''_{xxx}   \right) . $$
We start observing that, since $\Psi_F $ is a negative decreasing function, 
$$64 \Psi_F^4(x) \geqslant 64 \Psi_F^4(0)= 64 (4 \ln2)^4
 = 3782.\cdots.$$

\noindent
Since $(\Psi_F)'_x(x)$ is a not positive and decreasing function (point \ref{proprietà.derivata.1W}),
then $192 \Psi^2_F  (\Psi_F)'_x$ enjoys the same property. In particular
$$192 \Psi^2_F  (\Psi_F)'_x \geqslant 192 \Psi^2_F(x_{max})  (\Psi_F)'_x(x_{max})= -282.\cdots.$$

\noindent
From the previous observations it follows that $64  \Psi_F (\Psi_F)''_{xx} \geqslant 0,$
$48 ((\Psi_F)'_x)^2 \geqslant 0$ and they are increasing functions. 

\noindent
As for the last summand which appears in the expression of  $W'''_{xxx},$ we observe that
$$  (\Psi_F)'''_{xxx}= 2 \sum_{k=1}^\infty \frac{1}{(2k-1)!}\left(\frac{1}{2^{(2k+2)}}
 \psi^{(2k+2)}\left(\frac 1 2 \right)  - \psi^{(2k+2)}(1)\right) x^{2k-1}.$$ It is a not 
positive and decreasing function.
So we can write 
$$ 8 (\Psi_F)'''_{xxx}(x)\geqslant 8 (\Psi_F)'''_{xxx}(x_{max})=-19.9\cdots.$$

\noindent
We can conclude  $W'''_{xxx}(x)>0.$
Furthermore from our observations it follows that
$$W'''_{xxx}(x)\leqslant 
\left(  16 \Psi_F^4(x_{max})+24 ((\Psi_F)'_x )^2 (x_{max})\right.$$
$$\left.+ 16  \Psi_F(x_{max})  
 (\Psi_F)''_{xx}(x_{max})    \right)  <C_W$$ with $C_W=1125.$\\

\item  It holds that $Y'_x(x)=W(x)+x\,W'_x(x) .$ 
From the points \ref{proprietà.RF2}, \ref{proprietà.derivata.1W} and \ref{proprietà.derivata.2W}
we know that  $W(x)$ is a negative increasing function 
and $W'_x(x)$ is positive and decreasing for $x \in [0,x_{max}] .$
So we can write $W(x)\leqslant W(x_{max})= -4.1\cdots$
and  $ W'_x(x)\leqslant W'_x(0)=64 \ln^2 2 = 30.7\cdots.$
Then $Y'_x(x) \leqslant% \max W(x)+x_{max} \max W'_x(x)=
W(x_{max})+x_{max} W'_x(0)<0.$\\

\item It holds that  $Y''_{xx}(x)=2W'_x(x)+x\,W''_{xx}(x).$ 
From the points \ref{proprietà.derivata.1W}, \ref{proprietà.derivata.2W}
and \ref{proprietà.derivata.3W}
we know that  $W'_x(x)$ is a positive decreasing function 
and $W''_{xx}(x)$ is negative and increasing.
So we can write $W'_x(x)\geqslant W'_x(x_{max})= 22.\cdots.$
and  $W''_{xx}(x)\geqslant W''_{xx}(0)=-64 \ln^3 4 -6 \zeta(3) = -177.\cdots.$
Then $Y'_x(x) \geqslant %2 \min W'_x (x)+x_{max} \min W''_{xx}(x)=
2W'_x(x_{max})+x_{max} W''_{xx}(0)>0.$\\

\item It holds that  $Y'''_{xxx}(x)=3W''_{xx}(x)+x\,W'''_{xxx}(x).$ 
From the points  \ref{proprietà.derivata.2W}
and \ref{proprietà.derivata.3W} we know that  $W''_{xx}(x)$ is a negative increasing function 
and $0<W'''_{xxx}(x)<C_W.$ Then
$Y'''_{xxx}(x)\leqslant 3 W''_{xx}(x_{max})+ x_{max} C_W<0.$

\end{enumerate}
\hfill \qed

\noindent
\begin{proposition}
\label{minorazione.secondo.membro}
For all the  values of $l,x,y,z$ for which $T_2(l,x,y,z)$ is real, 
it holds that 
$$T_2(l,x,y,z) \leqslant \frac {1 + Cl^2 x}{l-1}, $$ where $C=-4 \ln 2.$
\end{proposition}
\noindent
{\bf Proof.} 
The epression of $T_2$ is given by (\ref{secondo.addendo}). We rewrite it in the following way
$$T_2=
%Q(x,y,z,l)=
%\sqrt{\left( \frac{l}{l-1}\right)^2 F^2 \left( x\right)  I^2\left( z\right)- 
%\frac{l+1}{l-1} L ( y)  G (z ) }=$$
%$$
\frac{1}{l-1}\sqrt{ l^2  F^2 \left( x\right)  I^2\left( z\right)-
 (l^2-1) L ( y)  I (z ) }.$$ 
We start studying the case of $T_2$ non zero.
If $1+ \bar R(x,y,z,l)$ is the Mac-Laurin series of the function under the square root then
we can write
$$T_2=\frac{1}{l-1}\sqrt{1+ \bar R(x,y,z,l)},$$
where $\bar R(x,y,z,l)=l^2 (R_{F^2}(c_1x)x (1+R_{I^2}(c_2 z)z) + R_{I^2}(c_2 z)z)-
(l^2-1)R_{LI}(y,z),$ and $c_1,c_2 \in (0,1).$ 
Thanks to the points \ref{proprietà.RLI}, \ref{proprietà.RI2}, \ref{proprietà.RF2}
and  \ref{proprietà.derivata.1W} of proposition \ref{proprietà}, we know 
that $R_{LI}(y,z)\geqslant 0,$  
%(point \ref{proprietà.RLI} of proposition \ref{proprietà}) 
 $R_{I^2}(x) \leqslant 0$  %(point \ref{proprietà.RI2} of proposition \ref{proprietà}) 
and that $R_{F^2}(x) $ is a 
negative increasing function, 
%(points \ref{proprietà.RF2} and  \ref{proprietà.derivata.1W}
%of proposition \ref{proprietà}), 
so $R_{F^2}(c_1 x) \leqslant R_{F^2}(x).$  
We can conclude that,
if we set $$R(x,z,l)=l^2 R_{F^2}(x) x (1+R_{I^2}(c_2 z)z),$$ 
then
$$\bar R(x,y,z,l) \leqslant l^2 R_{F^2}(c_1 x) x (1+R_{I^2}(c_2 z)z) \leqslant R(x,z,l), $$ 
%$$-Q(x,y,z,l)
$$T_2= \frac{1}{l-1}\sqrt{1+ \bar R(x,y,z,l)} \leqslant \frac{1}{l-1} \sqrt{1+ R(x,z,l)}.$$

\noindent
We know that 
$$\sqrt{1+f(x)}=\sqrt{1+f(0)}+ \frac {f'_t(t)}{2 \sqrt{1+f(t)}}|_{t=cx}x,$$
where $c \in (0,1).$ 
If  we apply this result to the function $f(x)= R(x,z,l),$ we get 

$$T_2 \leqslant \frac {\sqrt{1+  R (x,z,l)}}{l-1}=
\frac{1}{l-1}\left( \sqrt {1+R (0,z,l)} +
 \frac{R'_t(t,z,l) }{2\sqrt {1+R(t,z,l)} }|_{t=cx}  x\right),$$
where $c \in (0,1).$ 
We observe that $R(0,z,l)=0.$ %l^2 R_{I^2}(z) z^2-(l^2-1)R_{LI}(y,z)\leqslant 0.$
Then $$T_2 \leqslant \frac{1}{l-1}\left( 1+ \frac{R'_t(t,z,l) }
{2\sqrt {1+R(t,z,l)} }|_{t=cx}  x\right).$$
%$$=\frac{1}{l-1}\left( -1- \frac{R'_t(t,y,z,l) }{\sqrt {1+R(t,y,z,l)} }|_{t=cx}  x\right)$$
The proof will be completed after having proved the following result.
\hfill \qed

\begin{proposition}
\label{stima.minorazione.secondo.membro}
Under the same hypotheses of proposition \ref{minorazione.secondo.membro}
%For all the admissible values of $t,y,z,$
$$\frac{R'_t(t,z,l) }{2\sqrt {1+R(t,z,l)} }\leqslant Cl^2,$$
where $C=-4 \ln 2.$ 
\end{proposition}

\noindent
{\bf Proof.} We set $H(z,l)=l^2 (1+R_{I^2}(c_2 z)z)\leqslant l^2$ 
and $Y(t)=R_{F^2}(t) t.$
From the expression of $R(t,z,l)=H(z,l) Y(t),$ it follows that 
$R'_t(t,z,l)= H(z,l)Y_t'(t).$ 
Furthermore we can write 
$$ \frac{R'_t(t,z,l) }{2\sqrt {1+R(t,z,l)} }=\frac{H(z,l) Y_t'(t) }{2\sqrt {1+H(z,l) Y(t)} }.$$
We know from proposition \ref{proprietà} that $Y(t)\leqslant 0$ and $Y'_t(t)<0,$ 
then $R'_t(t,z,l)=H(z,l) Y'_t(t) \geqslant l^2 Y'_t(t),$
and $$-\frac{1}{2\sqrt{1+R(t,z,l)}}\geqslant -\frac{1}{2\sqrt{1+l^2 Y(t)}}.$$
We can conclude that 
$$-\frac{R'_t(t,z,l) }{2\sqrt {1+R(t,z,l)} }\geqslant  -\frac{l^2 Y'_t(t) }{2\sqrt {1+l^2Y(t)} } .$$
We shall show that the function on the right side is increasing with 
respect to the variable $t.$
The derivative with respect to the variable $t$ of this function 
is 
$$D(t,l)= -\frac{l^2}{2} \frac
{Y_{tt}''\sqrt{1+l^2 Y } - l^2 \frac{ (Y_t')^2 }
 {2\sqrt {1+l^2 Y} } } {1+l^2 Y }.$$
\noindent
We want to study the sign of $D(t,l).$ We start observing that 
 $1+l^2 Y \geqslant 1+\bar R >0.$
So it is sufficient to prove that the quantity 
$$ E(t,l)=2 Y_{tt}'' (1+l^2 Y) -l^2 (Y_t')^2 $$ is always not positive.
It holds that 
$$ Y'_t(t)= R_{F^2}(t)+t  (R_{F^2})_t'(t)$$ and 
$$ Y''_{tt}(t)= 2 (R_{F^2})'_t (t)+t  (R_{F^2})_{tt}''(t). $$
Then
$ Y(0)=0,$
 $Y'_t(0)=R_{F^2}(0)=2C$ and  
$Y''_{tt}(0)=2 (R_{F^2})'_t (0)=8 \Psi_F(0)^2=8 C^2.$
Furthermore we observe that $l\geqslant 2.$
So $$E(0,l)=16 C^2- 4l^2 C^2  \leqslant 0$$ and the equality holds if $l=2.$
The next step is to show that $E'_t(t,l)\leqslant 0.$
It is possible to find the following equality
$$ E'_t(t,l) =Y'''_{ttt}(1+l^2 Y) $$
Observing that $1+l^2 Y >0$ and $Y'''_{ttt}< 0$ (see the point \ref{proprietà.derivata.3Y}
of proposition \ref{proprietà}), we can conclude that $D(t,l)\geqslant 0$
(the equality holding if $t=0$). We have showed that
$$ -\frac{l^2 Y'_t(t) }{2 \sqrt {1+l^2Y(t)} }$$ is  a non decreasing function.
It gets the minimum for $t=0$ and its value is $-C l^2.$
Then 
$$-\frac{R'_t(t,z,l) }{2\sqrt {1+R(t,z,l)} }\geqslant -Cl^2,$$
and the proof of proposition \ref{stima.minorazione.secondo.membro}
 is completed. \hfill \qed \\

\noindent
To achieve the proof of proposition \ref{minorazione.secondo.membro},
we need  show that the statement continues to hold also for 
values of $l,x,y,z$ for which $T_2=0.$ 
To get this aim it is sufficient to observe that we can extend 
the result obtained under the hypothesis $T_2 >0$ 
for continuity.\\ 

\noindent
As for the first summand which appears in the expression of $t^2_-,$ that is 
$T_1,$ the following result holds. 
\begin{proposition}
\label{minorazione.primo.membro}
For all the admissible values of $x,z,$ it holds that
$$T_1\geqslant \frac l {l-1} (1+C x)$$
where $C=-4 \ln 2.$
 \end{proposition}
 \noindent
{\bf Proof.}
We recall that
$$T_1=\frac l {l-1} F(x) I(z)=\frac l {l-1}( 1 + R_{FI}(x,z)).$$ 
%where $R_{FI}(x,z) $= R_F(d_1 x) x +  R_I(d_2 y) y + R_F(d_1 x)  R_I(d_2 y) x y,$
%where $d_1,d_2 \in (0,1).$
Thanks to the point \ref{proprietà.RFI} of proposition \ref{proprietà} we have  
$R_{FI}(x,z) \geqslant C x.$
%Since  $R_FG(x,z)$ is an increasing function with respect to the variables $x$ and $z,$ 
% is a not positive and decreasing function and $R_F(t)$ is a 
%negative and increasing function, then 
%$R_{FI}(x,z) \geqslant   R_F(d_1 x) x + R_I(d_2 z) z \geqslant R_F(0) x=C x.$
%$R_{FI}(x,z) \geqslant R_F(0) x=C x.$ 
Then the result is immediate. \hfill \qed\\

\noindent
The following result gives the estimate of $t_-^2.$

\begin{proposition}
\label{minorazione.totale}
For all the  values of $x,y,z$ for which $t^2_- \in \R,$ it holds
$$t^2_- \geqslant 1-C l x,$$
where $C=-4 \ln 2.$
 \end{proposition}

\noindent 
We recall that
$t^2_-= T_1-T_2 .$ Thanks to propositions \ref{minorazione.secondo.membro} 
and \ref{minorazione.primo.membro} we get
$$t^2_-\geqslant \frac l {l-1} (1+C x)+ \frac 1 {l-1}(-1-Cl^2 x)=$$
$$ 1+ \left( \frac{C l}{l-1}- \frac{ C l^2}{l-1}\right) x =1+
\left (\frac{-Cl}{l-1}(l-1) \right) x =1-C l x.$$ \hfill \qed

\noindent
Now we turn our attention to the function $t_3.$ 
We recall that $s_{max}=\frac 1 {38}.$
\begin{proposition}
\label{stima.t3}
For  $s \in [0, s_{max}]$ 
  $$t^2_3(s)\leqslant 1+\frac{7}{2} s.$$
\end{proposition}
%\noindent
%As for $t_3^2,$ its Mac-Laurin series it is given by
%$$t^2_3(s)=1-\frac{1}{2}C s+D(ms)s^2,$$ where $m\in (0,1)$
\noindent
{\bf Proof.} 
We recall that 
$$t_3^2(s)=T(s)=\frac{4}{(1-s)^2} \left( \frac{\G(1+s)}{\G(1-s)}\right)^ {3}
 \frac{\G(1-2s)}{\G(1+2s)} \,\left( \frac{\G(3/2-s/2)}{\G(1/2+s/2)}\right)^ {2}.$$

\noindent
It holds that
$$T'_s(s)=\frac{1}{(1-s)} T(s) B(s),$$
%$$T'_s(s)=\frac{4}{(1-s)^3} \left( \frac{\G(1+s)}{\G(1-s)}\right)^ {3}
%\frac{\G(1-2s)}{\G(1+2s)} \,\left( %\frac{\G(3/2-s/2)}{\G(1/2+s/2)}\right)^ {2} B(s),$$
\noindent
where 
$$B(s)=2+(1-s) \left( -2\psi(1-2s)-2\psi(1+2s) +3\psi(1-s)+3\psi(1+s)- \right.$$
$$\left. -\psi \left(\frac 3 2- \frac s 2\right)- \psi \left(\frac 1 2+ \frac s 2\right)\right) .$$%Prima era 1/2-s/2.
To complete the proof we need the following result.
%\noindent
%We set $C(s)=-2\psi(1-2s)-2\psi(1+2s) +3\psi(1-s)+3\psi(1+s)- \psi %\left(\frac 3 2-
%\frac s 2\right)- \psi \left(\frac 1 2- \frac s 2\right).$

\begin{proposition}
If  $s\in [0,s_{max}]$ then $1<B(s)<3.$\\  
\end{proposition}

\noindent
{\bf Proof.} 
We observe that for $s \in [0,s_{max}]$
$$0<\psi \left(\frac 3 2- \frac s 2\right)<\psi\left(\frac 3 2 \right)=0.036\cdots,
 \quad
 \frac 3 2<- \psi \left(\frac 1 2+ \frac s 2\right)<
 -\psi\left(\frac 1 2 \right) < 2.$$
We can conclude that 
 $$1<- \psi \left(\frac 1 2+ \frac s 2\right)
 - \psi \left(\frac 3 2- \frac s 2\right)<2.$$
%$$ -\psi\left(\frac 1 2 \right) <
%- \psi \left(\frac 1 2- \frac s 2\right)< 2.$$
Furthermore
$$\psi(1-s)+\psi(1+s)=2 \sum_{k \geqslant 0} \frac{\psi^{(2k)}(1)}{(2k)!}s^{2k}, $$
from which it follows that 
$$D(s)= -2\psi(1-2s)-2\psi(1+2s) +3\psi(1-s)+3\psi(1+s)=2 \sum_{k \geqslant 0} 
\frac{\psi^{(2k)}(1)}{(2k)!}s^{2k}(3-2^{2k+1}).$$
If $k \geqslant 1$ then  $3-2^{2k+1} <0$ and $\psi^{(2k)}(1)<0$ 
(see formula 6.4.2 of \cite{AS})  then 
$$ 2 \psi(1)=-2 \g_{EM}= D(0) \leqslant  D(s) \leqslant  D(s_{max})=-1.146\cdots,$$
%with $D(s_{max})>0.$
%$$ -2\psi(1-2s)-2\psi(1+2s) +3\psi(1-s)+3\psi(1+s)\leqslant$$
%$$\leqslant    \psi(1-s)+\psi(1+s)=
%2 \sum_{k \geqslant 0} \frac{\psi^{(2k)}(1)}{(2k)!}s^{2k} \leqslant 2 \psi(1).$$
where $\g_{EM}=0.577\cdots$ is the Euler-Mascheroni constant.
\noindent
%So  $0< -\psi\left(\frac 3 2 \right)-\psi\left(\frac 1 2 \right)+ D(0)< %C(s) \leqslant D(s_{max})+2$ and  $2< B(s)\leqslant 2 + D(s_{max})< 2.85.$ 
So  $$ 1<B(s)\leqslant 2 +(1-s)(2+ D(s_{max}))< 4+D(s_{max})<3.$$
 \qed \\

\noindent
Since $B(s)>0$ then $T(s)$ is an increasing function and  we can deduce 
that $$T'(s) = \frac{1}{1-s}T(s)B(s)\leqslant \frac{3}{1-s_{max}}T(s_{max})  <7/2.$$ 

\noindent
The Mac-Laurin series of  order zero of $T(s) $ is  $1+T'_s(cs) s,$ where $c\in (0,1).$  
So it is immediate to conclude that 
$$T(s)\leqslant 1+\frac{7}{2} s.$$  \qed

%\noindent
%We want to remark that with more work it is possible
%to show that $$T(s)\leqslant 1-C s.$$
\noindent
The following proposition shows that 
the eventual solutions $t_{+}(l,g) \geqslant t_{-}(l,g)$ of the equation 
(\ref{equationt}) are always bigger than $t_3.$

\begin{proposition}
$t_3(\frac l {g+1})<t_-(l,g)$ for $g \geqslant 1$ and   
$2\leqslant l \leqslant g-1$ such that $t_-(l,g) \in \R.$
\end{proposition}

\noindent
{\bf Proof.} 
From our observations, it is sufficient to show that  $t^2_3(s)<t^2_-(l,g)$ holds
for $g \geqslant 38.$
Propositions \ref{minorazione.totale} and \ref{stima.t3} assure that 
$$t^2_- \geqslant 1- C lx,$$
$$t^2_3(s)\leqslant 1+\frac{7}{2} s.$$
We recall that $x=ls$ and $2 \leqslant l \leqslant g-1.$ Then the result is obvious.
\hfill \qed

\section{The index and the nullity of the Costa-Hoffman-Meeks surfaces}
\label{statements}

We start recalling some results described in previous sections.
We denoted by $G_t,$ $t\in (0,+\infty),$ a deformation of the map 
$G$ defined by (\ref{definizione.G}). 
Thanks to  theorem \ref{nullity.37},  ${\rm Nul}(G_t)>3$ only if $t$ 
assumes special values.
If $2\leqslant g \leqslant 37$ these values are $t_1,t_2,t_3.$
If $g\geqslant 38$ there are additional values.
They are the positive solutions of the  equation (\ref{equationt}).
We denoted them by $t_\pm(l,g),$ where $2 \leqslant l \leqslant g-1,$ 
and for definition $ t_+\geqslant t_-.$
In  previous section we have proved that the inequality
$t_3 (s)< t_-(l,g) $ holds.
S. Nayatani  showed in \cite{N1} that  $t_3>t_2$ for   $g \geqslant 2.$
%and theorem \ref{nullity.37}.
We can conclude that  no one of the $t_\pm$
can be equal to $t_2.$ As consequence ${\rm Nul}(G_{t_2})$ continues
to be equal to 4 also for $g \geqslant 38,$ because ${\rm dim}\,H(G_{t_2})$
is equal to 1 for all $g \geqslant 2.$\\

\noindent
We recall that $M_g$ denotes the Costa-Hoffman-Meeks surface of genus $g.$
Since the extended Gauss map of  $M_g$ is exactly $G_{t_2},$ and
taking into account  the result of S. Nayatani about the Costa surface 
 (theorem \ref{nullindex.Costa}) showed in \cite{N2}, we have proved 
 theorem \ref{teorema.nullity}.\\
 
 %the following result. OLD
%\begin{theorem}
%\label{teorema.nullity}
%The null space of the Jacobi operator of $M_g$
%%Costa-Hoffman-Meeks surface 
%has dimension equal to $4$ for all $g \geqslant 1.$
%\end{theorem}
%
%\noindent
%Using the definition of non degeneracy given in \cite{PR}, we can also rephrase 
%this result giving the following statement. 
%
%\begin{corollary}
%\label{corollario.nondeg}
%The surface $M_g$
%%Costa-Hoffman-Meeks surface 
% is non degenerate for all $g \geqslant 1.$
%\end{corollary}

\noindent
Now we turn our attention to the results relative to the
index of the map $G_t.$ We recall that $\Sigma_g$ denotes the compactification
of $M_g.$ 
%We recall that the index of a complete oriented minimal surface is 
%the index of its extended Gauss map $G.$ 
%We have denoted it with ${\rm Ind}(G).$
%It is defined to be the number of negative eigenvalues of 
%the operator $L_G.$
S. Nayatani proved in \cite{N1} the following result.

\begin{theorem}
\label{teorema.indice}
Let $G:\Sigma_g \to S^2$ be the holomorphic map defined by (\ref{definizione.G}). 
If $2 \leqslant g \leqslant 37,$ then
$$
{\rm Ind}(G_t)= \left\lbrace 
\begin{array}{cl}
2g+3 & \mbox{if } \,\,t \leqslant t_1,\, t_2 \leqslant t < t_3,\, t > t_3,\\
2g+4 & \mbox{if } \,\,t_1 < t < t_2, \\
2g+2 & \mbox{if } \,\,t=t_3.
\end{array}
\right.
$$
\end{theorem}

\noindent
For $t=t_1,t_2,t_3$ we have ${\rm Nul}(G_t)>3,$ that is 
the kernel of $L_{G_t}$ contains at least one non linear element.
The eigenvalue associated  to this function is zero.
%So the a null eigenvalues of $L_{G_t}.$
The proof of  theorem
\ref{teorema.indice} is based on the analysis of the behaviour of 
these null eigenvalues under a variation of the value of $t.$
Let's suppose that $t \neq t_1,t_2,t_3$ but remaining in a neighbourhood 
of one of these values. For example we choose $t_1.$
Then the eigenvalue $E$ that before the variation was associated to a 
non linear element of $N(G_{t_1}),$ is no more equal to zero. 
To compute the index, it was necessary to understand which
is the sign assumed by $E,$ respectively  for $t>t_1$ 
and $t<t_1.$ Similar considerations are applicable to 
the eigenvalues associated with $t_2$ and $t_3.$
See \cite{N1} for the details.\\

\noindent
If $g \geqslant 38,$ we have just proved that the other 
values for which ${\rm Nul}(G_t)>3$ are bigger than $t_3.$
The presence of these additional values $t_\pm$ does not 
influence the value of ${\rm Ind}(G_t)$ if $t \leqslant t_3.$
In other terms theorem \ref{teorema.indice}
continues to hold for $g\geqslant 38$ if we consider $0< t \leqslant t_3.$
Taking into account also the result of S. Nayatani about the Costa 
surface  ($g=1$) showed in 
\cite{N2}, we have proved theorem \ref{indice.costa}.

%OLD
%we can give the following statement
%
%\begin{theorem}
%\label{indice.costa}
%For all $g\geqslant 1$ the index of the Gauss map of $M_g$ is equal to $2g+3.$    
%\end{theorem}

\section{Appendix}
\label{appendix}
This section contains some additional details of the 
computations made by S. Nayatani.
%A parallel reading of this section
%with sections \ref{caso.minore38} and \ref{caso.38}
%is strongly recommended.

\subsection{Divisors and  Riemann-Roch theorem}
\label{RiemannRoch}
Here we  introduce some definitions and concepts
 of the algebraic geometry. See for example \cite{DS}.\\

\noindent
Let $\Sigma_g$ be a compact Riemann surface of genus $g.$ 
A divisor on $\Sigma_g$ is a finite formal sum
of integer multiples of points of $\Sigma_g,$
$$D=\sum_{x \in \Sigma_g} n_x x, \quad n_x \in \Z, n_x=0 \quad \hbox{{\rm for almost all}}\quad x.$$
The set of the divisors on $\Sigma_g$ is denoted by Div$( \Sigma_g).$
The degree of a divisor is the integer deg$(D)=\sum n_x.$ \\

\noindent
Let $	{\mathbb C}(\Sigma_g)$ be the field of the meromorphic functions on $\Sigma_g$ 
and let ${\mathbb C}(\Sigma_g)^*$ be its multiplicative group of nonzero elements.
Every $f \in {\mathbb C}(\Sigma_g)^*$ has a divisor 
$$ {\rm div}(f)=\sum \nu_x(f) x, $$  
where $ \nu_x(f)$ denotes the order of $f$ at $x.$\\

\noindent
Let $\omega$ be a nonzero meromorphic differential $n$-form on $\Sigma_g.$
Then $\omega$ has a local representation $\omega_x= f_x(z)(dz)^n $ about 
each point $x$ of $\Sigma_g,$ where $z$ is the local coordinate about $x$
and $f_x(z) \in {\mathbb C}(\Sigma_g)^*.$ 
So we can define in a natural way $\nu_x(\omega)=\nu_0(f_x)$ and also
associate  a divisor with a differential form: 
$$\hbox{{\rm div}}(\omega)=\sum \nu_x(\omega) x.$$
A canonical divisor on $\Sigma_g$ is a divisor of the form div$(\omega)$ where
$\omega$ is a nonzero meromorphic differential form.  \\

\noindent
Let $D \in {\rm div}(\Sigma_g).$
We denote by $H^{0,n}(D)$  the vector space of the meromorphic
differential $n$-forms $\omega$ 
such that $${\rm div}(\omega)+D \geqslant 0.$$ 
In other terms, if 
$D={\rm div}(\eta),$ with $\eta$ differential form with
local representation $\eta_x= g_x(z)(dz)^n ,$ then 
the elements of $H^{0,n}(D)$ 
are the differential forms
$\omega$ having a local representation $\omega_x= f_x(z)(dz)^n $ with
  $f_x \in {\mathbb C}(\Sigma_g)$ vanishing to high enough order
to make the product $f  g$ holomorphic. 
We set ${\rm dim}\,H^{0,n}(D)=\ell(D).$ \\

\noindent
We are ready to state the following result. 

\begin{theorem}[Riemann-Roch]
\label{riemann.roch}
Let $\Sigma_g$ be a compact Riemann surface of genus $g.$ Let $k_{\Sigma_g}$ be a canonical divisor
on $\Sigma.$ Then for any divisor $D \in {\rm Div}(\Sigma_g),$
$$\ell(D)={\rm deg}(D)-g+1+ \ell(k_{\Sigma_g}-D).$$   
\end{theorem}

\noindent
The next result gives information about the canonical divisor 
and a simpler version of Riemann-Roch theorem for 
divisors of large enough order.
\begin{corollary}
\label{corollario.riemann.roch}
Let $\Sigma_g,g,D,k_{\Sigma_g}$ as above. 
\begin{itemize}
\item ${\rm deg}(k_{\Sigma_g})=2g-2,$
\item If ${\rm deg}(D)>2g-2$ then $\ell(k_{\Sigma_g}-D)=0.$
Equivalently
$\ell(D)={\rm deg}(D)-g+1.$
\end{itemize}
\end{corollary}

\subsection{The determination of a basis of differential forms 
with null residue at the ramification points}
\label{base}
The ramification points (or branch points) of $\gamma(w)=w$ are 
 the zeroes of 
$$\frac{dw}{dz}=\frac{g+2}{g+1}\,
\frac{z^{g-1}(z^2-A^2)}{w^{g}}=
\frac{g+2}{g+1}\, \frac{z^{g-1}(z^2-A^2)}{(z^{g}(z^2-1))^{\frac{g}{g+1}}},$$
with $A=\sqrt{\frac{g}{g+2}},$ where $g$ denotes the genus, the pole of 
$\gamma$ and the origin of $\C^2.$  
That is $Q_0=(0,0),$ $P_\infty=(\infty,\infty),$ 
$P_m=(A,B_m)$ and $S_m=(-A,C_m)$ for $m=0,\ldots,g,$ 
where $B_m,C_m$ denote, respectively, the $m$-th complex value  of 
$\sqrt[g+1]{A^g(A^2-1)}$ and $\sqrt[g+1]{(-A)^g(A^2-1)}.$ 
We  have set $P_\pm=(\pm 1,0).$
We recall that 
\begin{equation}
\hat H(G)= \left\lbrace  \sigma \in H^0(2k_{\Sigma_g}+R(G)) \,| \, Res_{r_i} \frac {\sigma} {dw}=0, i=1,\ldots, \mu\right\rbrace,
\end{equation}
where $k_{\Sigma_g}$ is a canonical divisor of ${\Sigma_g}$
and $R(G) =\sum_1^\mu r_i$ is the ramification divisor of $G.$
% where $r_i$ are the branch points of $\g(w)=w.$
In our case it is given by $R(G)=Q_0+P_\infty+\sum_{m=0}^{g}(P_m+S_m).$
Furthermore it holds $\hat H(G)=\hat H(G_t).$\\

%$$H(G_t)=\left\lbrace  \sigma \in \hat H(G) \, |\,Re 
%\int_\a (1-t^2 g^2,i(1+t^2 g^2),2tg)
%\frac {\sigma} {dw}=0, \, \forall \a \in H_1(\Sigma,\Z) \right\rbrace,$$ 

\noindent
As for the canonical divisor $k_{\Sigma_g},$ we consider 
$k_{\Sigma_g}=
(g-1)P_++(g-1)P_-.$
We observe that ${\rm deg}(k_{\Sigma_g})=2g-2$ like
stated by corollary \ref{corollario.riemann.roch}.\\

\noindent
To study the space $\hat{H}(G_t)$ we need understand
which are the elements of the space $H^0(2k_{\Sigma_g}+R(G)).$
Taking into account the definitions of $k_{\Sigma_g}$ and 
$R(G),$ then 
$2k_{\Sigma_g}+R(G)=2(g-1)P_+ +2(g-1)P_- + 
Q_0+P_\infty + \sum_{m=0}^g P_m + \sum_{m=0}^g S_m.$ 
%From the definition of $H^0,$ 
Among  
%We deduce that
the quadratic differentials $\sigma$ that are 
in $H^0(2k_{\Sigma_g}+R(G)),$ 
we consider the ones having one of the following forms:
%can have two possible structures:
\begin{equation}
\label{famiglia1}
z^kw^j\left( \frac {dz} w  \right) ^2,
\end{equation}
\begin{equation}
\label{famiglia2}
z^kw^j \frac 1 {z\pm A} \left( \frac {dz} w  \right) ^2.
\end{equation}
In fact from the definition of $H^0,$ it follows that the 
quadratic differentials to consider can have
 a pole of order $0$ (differentials of  type 
(\ref{famiglia1})) or of order $1$ (differentials of  type 
(\ref{famiglia2})) at $P_m$ and $S_m$ 
for $k=0,\ldots,g.$
%Furthermore we observe that the quadratic differential 
%$\left( \frac {dz} w  \right) ^2$ has 
%\begin{itemize}
%\item a zero of order $1$ at $Q_0$
%\item a zero of order $1$ at $P_\infty$
%\item a zero of order $2(g-1)$ at $P_+$ and $P_-.$
%\end{itemize}
%The quadratic differential 
%$\frac{1}{z \pm A}\left( \frac {dz} w  \right) ^2$ 
%has 
%\begin{itemize}
%\item a zero of order $1$ at $Q_0$
%\item a zero of order $2(g-1)$ at $P_+$ and $P_-$
%\item a zero of order $g+2$ at $P_\infty$
%\item 
%a pole of order $1$ at $P_m$ if we consider sign "-" or at $S_m$
%if we consider the sign "+".
%\end{itemize}
We will  determine separately which are the 
differential forms of type (\ref{famiglia1}) and
(\ref{famiglia2}) belonging to $H^0(2k_{\Sigma_g}+R(G))$. 
To select the differential forms 
of type (\ref{famiglia2}) it is convenient 
to introduce an auxiliar divisor. 
$$D= Q_0+ (g+2) P_\infty + 2(g-1)P_+ +2(g-1)P_- .$$
Actually to determine the differential forms of type (\ref{famiglia2})
which belong to $H^0(2k_{\Sigma_g}+R(G))$ is equivalent
to look for the differential forms of type  (\ref{famiglia1})
which are in $H^0(D).$
We observe that the elements of the vector space $H^0(D)$
after the multiplication by the factor 
$z \pm A$ are elements of $H^0(2k_{\Sigma_g}+R(G)).$
It is necessary to remark that to obtain 
a basis of  $H^0(2k_{\Sigma_g}+R(G)),$
we will not take into account the differentials
of $H^0(2k_{\Sigma_g}+R(G))$ that can be 
constructed from an element of $H^0(D)$
as described above. Otherwise the number of the 
founded differential forms would exceed the dimension
of $H^0(2k_{\Sigma_g}+R(G)),$ that we can 
compute as follows.
We observe that ${\rm deg}(2k_{\Sigma_g}+R(G))=6g.$ 
Then thanks to corollary \ref{corollario.riemann.roch} 
 %of Riemann-Roch theorem 
we conclude that ${\rm dim}\,H^0(2k_{\Sigma_g}+R(G))=5g+1.$
% and ${\rm dim}H^0(D)=4g.$ 
So the basis we are looking for counts $5g+1$ elements.
%We observe that 
From the observations made above we can deduce that
among the forms of  type (\ref{famiglia1}), we will consider
 the ones which satisfy  the following conditions
%$$ z^k w^l \in H^0(D)$$ or equivalently
$$
\left\{
\begin{array}{l}
 k(g+1)+j g \geqslant -1,\\
 j\geqslant -2 (g-1),\\
-k(g+1) -j(g+2) \geqslant -1.
\end{array}
\right.$$
These relations assure that a differential form $\omega$
of type
$z^kw^j   \left( \frac {dz} w  \right) ^2,$
%are holomorphic, respectively, at the points $Q_0,$ 
%$P_{\pm}$ and $P_{\infty}.$
satisfies ${\rm div }(\omega)+2k_{\Sigma_g}+R(G) \geqslant 0.$
These differentials can be classified in three families.
Each family is characterized by particular values
of $j$ and $k.$ That is
\begin{enumerate}
\item $j=-g+1,\ldots,0,1 $ and $k=-j,$
\item $j=2-2g,\ldots,-g $ and $k=-j,$ 
\item $j=2-2g,\ldots,-g $ and $k=-j-1.$ 
\end{enumerate}

\noindent
As for the forms  of  type (\ref{famiglia2}) we shall consider 
only the ones which satisfy 
%$$ z^k w^l \in H^0(2k_{\Sigma_g}+R(G))$$ or equivalently
$$\left\{
\begin{array}{l}
 k(g+1)+j g \geqslant -1,\\
 j\geqslant -2 (g-1),\\
-k(g+1) -j(g+2) \geqslant -(g+2).
\end{array}
\right.$$
These relations assure that a differential form
of type
$z^kw^j  \frac 1 {z\pm A} \left( \frac {dz} w  \right) ^2,$
satisfies ${\rm div }(\omega)+D \geqslant 0.$
%are holomorphic, respectively, at the points $Q_0,$ 
%$P_{\pm}$ and $P_{\infty}.$
We obtain that $j=-g+1,\ldots,0,1$ and $k=-j+1.$\\

\noindent
Since we are looking for a basis of a vector space
we can replace each couple of differentials
$\frac{f}{z-A}\left(\frac{dz}{w} \right)^2,
\frac{f}{z+A}\left(\frac{dz}{w} \right)^2$
by an appropriate linear combination.
We observe that 
$$\frac{1}{z-A}\pm \frac{1}{z+A}=
\left\lbrace \begin{array}{c}
\eta_1=\frac{z}{z^2-A^2}\\
\eta_2=\frac{1}{z^2-A^2}.
\end{array}\right.$$ 
So in the following we will work with 
the forms $f \eta_1 \left( \frac{dz}{w}\right)^2$ and  
$f \eta_2 \left( \frac{dz}{w}\right)^2,$
where $f=z^k w^j$ as described above.\\

\noindent
The $5g+1$  quadratic differentials we have found forms
a basis of $\hat{H}(G_t).$ The last step is to divide each elements 
of this basis by $dw.$  
After simple algebraic manipulations, we obtain the following $5g+1$ 
differential $1$-forms:
\begin{equation}
%\begin{enumerate}
\label{lista.candidati}
%\item
\begin{array}{ll}
 \frac{w^k}{z^{k-1}}\frac{dz}{(z^2-A^2)^2} & {\rm for}\quad k=-1,0,\ldots,g-1,\\
%\item 
\frac{w^k}{z^k}\frac{dz}{(z^2-A^2)^2} &  {\rm for}\quad k=-1,0,\ldots,g-1,\\
%\item 
\frac{w^k}{z^{k+1}}\frac{dz}{(z^2-A^2)} &{\rm for}\quad k=-1,0,\ldots,g-1,\\
%\item 
\frac{z^k}{w^{k+1}}\frac{dz}{(z^2-A^2)} &{\rm for}\quad k=1,\ldots,g-1,\\
%\item
 \frac{z^{k-1}}{w^{k+1}}\frac{dz}{(z^2-A^2)} & {\rm for}\quad k=1,\ldots,g-1.
%\end{enumerate}
\end{array}
\end{equation}

\noindent
Now it is necessary to select the $1$-forms having residue equal to zero 
at the points $Q_0,$ $P_m$ and $S_m$ with $m=0,\ldots,g.$ 
Thanks to the properties of symmetry of
the surface it is sufficient  to verify
the null residue condition at the points $Q_0,$ 
$P_1=(A, e^{\frac{2\pi i}{g+1}} \sqrt[g+1]{A^g(A^2-1)}).$ 
%$S_1=(-A, e^{\frac{2 \pi i}{g+1}} \sqrt[g+1]{(-A)^g(A^2-1)}).$
In fact from the coordinates of the points $P_m$ and $S_m,$
 we can deduce that for each $Q \in \{P_m,S_m , m=0,\ldots,g \}$ 
 there exists $n \in \{0,\ldots,2g+1 \}$ such that 
 $Q=\lambda^{n}(P_1),$ where $\lambda$ is the conformal 
 diffeomorphism described in lemma \ref{lemma.simmetrie}.
So we can state that the residue of an arbitrary form $\omega$ at the point
$Q$ is related to the residue at $P_1$ by
$${\rm Res}_{Q} \omega={\rm Res}_{P_1} (\lambda^{n-1})^*\omega.$$
Applying this result to the differential
forms of the list (\ref{lista.candidati}) and using the
the definition (\ref{lambdakappa}) of $\lambda,$ it is easy to obtain  
 that ${\rm Res}_{Q} \omega$ is equal to ${\rm Res}_{P_1} \omega$ times a 
 power of $\pm \rho.$ So if ${\rm Res}_{P_1}\omega=0$ then 
${\rm Res}_{Q}\omega=0.$\\

\noindent
Thanks to algebraic manipulations inspired by the simpler cases 
where $g=2,3,$ it is possible to find $3g$ linear independent
differential forms satisfying the null residue condition.
They constitute  the wanted basis.

$$\omega_k^{(1)}=\frac{z^{k-1}}{w^{k}}\frac{dz} w \quad {\rm for} \quad k=1,\ldots,g-1,$$

$$\omega_k^{(2)}=
\frac{z^{k-1}((k-2)z^2-kA^2)}
       {w^{k}(z^2-A^2)^2}dz \quad
 {\rm for} \quad k=0,\ldots,g,$$ 

$$\omega_k^{(3)}=\frac{z^{k-1}((k-2)z^2-kA^2)}
       {w^{k+1}(z^2-A^2)^2}dz
\quad {\rm for} \quad k=0,\ldots,g-1.$$

\subsection{The equations equivalent to the condition of existence of a branched minimal surface.}
\label{sistemi}
Let $\omega_1$ and $\omega_2$ two meromorphic differential forms on ${\Sigma_g}.$
We write $\omega_1 \sim \omega_2$ if there exists a meromorphic function
$f$ on ${\Sigma_g}$ such that $\omega_2=\omega_1+df.$  It is possible to prove
that:

$$\omega_k^{(2)}
\sim -\frac{k(g+2)}{2(g+1)} \frac{z^{k-1}}{w^{k}}dz \quad 
{\rm for}\quad k=0,\ldots,g,$$

$$\omega_k^{(3)}
\sim -\frac{(g+2)(g+k+2)}{2(g+1)} \frac{z^{k-1}}{w^{k+1}}dz
\quad {\rm for}\quad  k=0,\ldots,g-1.$$
\noindent
Using these relations we get:

$$\int_{\tilde \b}  \omega_k^{(1)}=-2 i \sin \frac{(k+1) \pi}{g+1} K_k, \quad
\int_{\tilde \b}  \omega_k^{(2)}=-\frac{(g+2)k}{2(g+1)} 2 i \sin \frac{k \pi}{g+1} I_k,$$
$$\int_{\tilde \b}  \omega_k^{(3)}=\frac{(g+2)(g+2-k)}{2(g+1)} 2 i \sin \frac{(k+1)
 \pi}{g+1} K_k,$$
$$\int_{\tilde \b} \g \omega_k^{(1)}=2 i \sin \frac{k \pi}{g+1} I_k, \quad
\int_{\tilde \b} \g \omega_k^{(2)}=\frac{(g+2)(g+2-k)}{2(g+1)} 2 i \sin \frac{(k-1)
 \pi}{g+1} J_k,$$
$$\int_{\tilde \b} \g \omega_k^{(3)}=-\frac{(g+2)k}{2(g+1)} 2 i \sin 
\frac{k \pi}{g+1} I_k,$$
$$\int_{\tilde \b} \g^2 \omega_k^{(1)}= 2 i \sin \frac{(k-1) \pi}{g+1} J_k, \quad
\int_{\tilde \b} \g^2 \omega_k^{(2)}=\frac{(g+2)(2g+4-k)}{2(g+1)} 2 i \sin 
\frac{(k-2) \pi}{g+1} L_k,$$
$$\int_{\tilde \b} \g^2 \omega_k^{(3)}=\frac{(g+2)(g+2-k)}{2(g+1)} 2 i 
\sin \frac{(k-1) \pi}{g+1} J_k.$$

\noindent
We recall that we must impose that $\omega =\sum_0^{g-1} c_k^{(1)} \omega_k^{(1)} +
\sum_0^{g} c_k^{(2)} \omega_k^{(2)}
+\sum_0^{g-1} c_k^{(3)} \omega_k^{(3)},$
where $c_k^{(i)} \in \C,$ satisfies
$$ \int_\a \omega= t^2  \overline{ \int_\a \g^2(w) \omega }, \quad
 Re \int_\a \g(w)  \omega =0$$
for $\a=\lambda^l\circ \tilde \beta$ for $l=0,\ldots,2g-1.$
Now it is convient to introduce some additional notation.\\

\noindent
Let
\begin{equation}
\label{matrice.rotazione}
{\mathcal L}=
\left[
\begin{array}{cc}
{\mathcal R}_\t & 0 \\
0&1
\end{array} \right]
\end{equation}
where ${\mathcal R}_\t$ is the rotation in the plane by $\t=g \pi/(g+1).$\\
%Then $$\lambda^*(\g \o_k^{(1)} )=(-1)^k \rho^{-k}\g \o_k^{(1)},$$
%$$\lambda^*(\g \o_k^{(2)} )=(-1)^k \rho^{-k+1}\g \o_k^{(2)},$$
%$$\lambda^*(\g \o_k^{(3)} )=(-1)^k \rho^{-k}\g \o_k^{(3)}.$$
If we denote $\Phi(\omega)=(1-\g^2,i(1+\g^2),2 \g)\omega,$
then it is possible to prove
$$\int_{{\lambda^l \circ \tilde \b}} \Phi(\omega)=
\int_{ \tilde \b} \lambda^* \Phi(\omega).$$
Since %we  want to apply this last relation to 
the differential form $\omega$ is linear combination 
of $\omega_k^{(j)},$ $j=1,2,3,$ 
 it is convenient to remark that:
$$\lambda^*\Phi( \omega_k^{(1)} )=(-1)^k \rho^{-k}{\mathcal L}\Phi( \omega_k^{(1)} ),$$
$$\lambda^*\Phi(\omega_k^{(2)} )=(-1)^k \rho^{-k+1}{\mathcal L}\Phi( \omega_k^{(2)} ),$$
$$\lambda^*\Phi(\omega_k^{(3)} )=(-1)^k \rho^{-k}{\mathcal L}\Phi( \omega_k^{(3)} ),$$
where $\rho=e^{i\frac{g\pi}{g+1}}.$
%We recall that 
% $$\o =\sum_0^{g-1} c_k^{(1)} \o_k^{(1)} +
%\sum_0^{g} c_k^{(2)} \o_k^{(2)}
%+\sum_0^{g-1} c_k^{(3)} \o_k^{(3)},$$
%where $c_k^{(i)} \in \C$ and $\rho=e^{i\frac{g\pi}{g+1}}.$ 
Then the equations 
$${\rm Re}\int_{\lambda^l \circ \tilde \b} (1-t^2\g^2,i(1+t^2 \g^2))\omega=0,\quad
{\rm for} \quad l=0,\ldots,2g-1,$$ are equivalent to:
$$Im\left[\sum_{k=0}^{g-1} \{(-1)^k \rho^{-k} \}^l f_k+ 
\sum_{k=1}^{g} \{(-1)^k \rho^{-(k-1)} \}^l p_k  \right]=$$
$$t^2 Im\left[\sum_{k=0,k\neq 1}^{g-1} \{(-1)^k \rho^{-k} \}^l h_k+ 
\sum_{k=0,k\neq 2}^{g} \{(-1)^k \rho^{-(k-1)} \}^l q_k  \right],$$
$$Re\left[\sum_{k=0}^{g-1} \{(-1)^k \rho^{-k} \}^l f_k+ 
\sum_{k=1}^{g} \{(-1)^k \rho^{-(k-1)} \}^l p_k  \right]=$$
$$-t^2 Re \left[\sum_{k=0,k\neq 1}^{g-1} \{(-1)^k \rho^{-k} \}^l h_k+ 
\sum_{k=0,k\neq 2}^{g} \{(-1)^k \rho^{-(k-1)} \}^l q_k  \right],$$
$l=0,\ldots,2g-1.$ These last equations can be arranged as in the systems
(\ref{sistema1}) and (\ref{sistema1.++}).
 The equations
$${\rm Re}\int_{\lambda^l \circ \tilde \b} 2t\g\omega=0,\quad
{\rm for} \quad l=0,\ldots,2g-1,$$
are equivalent to:
$$Im\left[\sum_{k=1}^{g-1} \{(-1)^k \rho^{-k} \}^l d_k+ 
\sum_{k=0,k\neq 1}^{g} \{(-1)^k \rho^{-(k-1)} \}^l e_k  \right]=0,$$
$l=0,\ldots,2g-1.$ 
These last equations can be arranged as in the systems
(\ref{sistema2}) and (\ref{sistema2.++}).

\end{document}